\documentclass[oneside, 10pt, a4paper]{article}
 \usepackage{amsmath}
 \usepackage{epsfig} 
  \usepackage{graphics} 
 
\usepackage{mathrsfs}
\usepackage{amsthm} 
\usepackage{amssymb,amsfonts}
\usepackage[english]{babel}
\usepackage[latin1]{inputenc}
\usepackage{amscd}
\usepackage{verbatim}
\newcommand{\famoustheo}[2][]{\mbox{}\par\textbf{#1.}\ \emph{#2}}
\DeclareMathOperator*{\Var}{Var}
\DeclareMathOperator*{\sign}{sign}
\providecommand{\pint}[1]{\ensuremath{\left[ #1 \right]}}
\providecommand{\abs}[1]{\ensuremath{\left\lvert #1 \right\rvert}}
\providecommand{\norm}[1]{\ensuremath{\left\Vert #1 \right\Vert}}
\providecommand{\mod}{mod}

\newtheorem{theorem}{Theorem}
\newtheorem{corollary}{Corollary}
\newtheorem{lemma}{Lemma}
\newtheorem{proposition}{Proposition}
\theoremstyle{definition}

\newtheorem{remark}{Remark}

\title{On the entropy of Japanese continued fractions}
\author{\small{\scshape{LAURA LUZZI, STEFANO MARMI}}}
\date{}

\begin{document}
\maketitle

\begin{abstract}
We consider a one-parameter family of expanding interval maps $\{T_{\alpha}\}_{\alpha \in [0,1]}$ (\emph{japanese continued fractions}) which include the Gauss map ($\alpha=1$) and the nearest integer and by-excess continued fraction maps ($\alpha=\frac{1}{2},\,\alpha=0$).
We prove that the Kolmogorov-Sinai entropy $h(\alpha)$ of these maps depends continuously on the parameter and that $h(\alpha) \to 0$ as $\alpha \to 0$.
Numerical results suggest that this convergence is not monotone and that the entropy function has infinitely many phase transitions and a self-similar structure.
Finally, we find the natural extension and the invariant densities of the maps $T_{\alpha}$ for $\alpha=\frac{1}{n}$.   
\end{abstract}

\section{Introduction}

Let $\alpha \in [0,1]$. We will consider the one-parameter family of maps $T_\alpha:I_{\alpha} \to I_{\alpha}$, where $I_{\alpha}=[\alpha-1,\alpha],$ defined by
\begin{equation*}
T_\alpha(x)=\abs{\frac{1}{x}}-\pint{\abs{\frac{1}{x}}+1-\alpha}
\end{equation*}

\begin{figure}[htb]
\begin{center}
\includegraphics[width=0.6\textwidth, height=0.6\textwidth, angle=270]{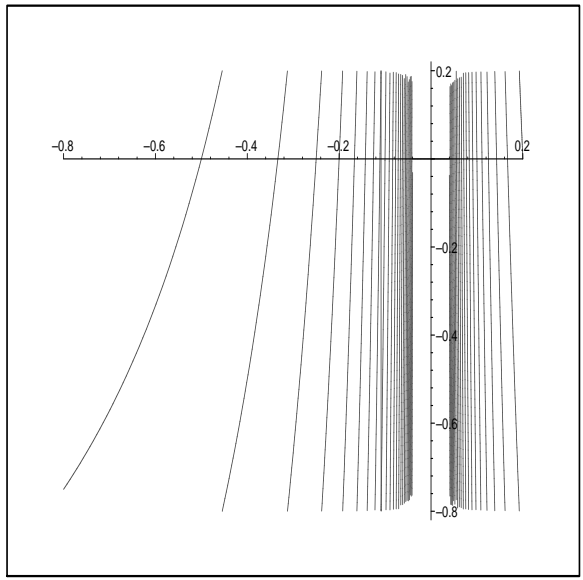}
\caption{Graph of the map $T_\alpha$ when $\alpha=0.2$}
\end{center}
\end{figure}

These dynamical systems generalize the Gauss map ($\alpha=1$) and the nearest integer continued fraction map ($\alpha=\frac{1}{2}$); they were introduced by H. Nakada \cite{Nak}. For all $\alpha \in (0,1]$ these maps are expanding and admit a unique absolutely continuous invariant probability measure $d\mu_{\alpha}=\rho_{\alpha}(x)dx$ (for a detailed proof in this particular case see for example \cite{BDV}). Nakada computed the invariant densities $\rho_{\alpha}$ for $\frac{1}{2} \leq \alpha \leq 1$ by finding an explicit representation of their natural extensions. The maps $\rho_{\alpha}$ are piecewise finite sums of linear fractional functions: 
\begin{multline}
\text{For } g < \alpha \leq 1, \\
\rho_{\alpha}(x)=\frac{1}{\log(1+\alpha)} \left(\chi_{\left[\alpha-1,\frac{1-\alpha}{\alpha}\right]}(x)\frac{1}{x+2}+\chi_{\left(\frac{1-\alpha}{\alpha},\alpha\right)}(x)\frac{1}{x+1} \right) \notag
\end{multline}
\begin{multline}
\text{For } \frac{1}{2} < \alpha \leq g, \\
\rho_{\alpha}(x)=\frac{1}{\log G} \left(\chi_{\left[\alpha-1,\frac{1-2\alpha}{\alpha}\right]}(x)\frac{1}{x+G+1}+ \right. \\
+\left. \chi_{\left(\frac{1-2\alpha}{\alpha},\frac{2\alpha-1}{1-\alpha}\right)}(x)\frac{1}{x+2} + \chi_{\left[\frac{2\alpha-1}{1-\alpha},\alpha\right)}(x)\frac{1}{x+G} \right) \notag
\end{multline}
where $g,G$ denote the golden numbers $\frac{\sqrt{5}-1}{2}$ and $\frac{\sqrt{5}+1}{2}$ respectively.\smallskip

\par The case $\sqrt{2}-1 \leq \alpha \leq \frac{1}{2}$ was later studied by Moussa, Cassa and Marmi \cite{MCM} for a slightly different version of the maps, that is  $M_\alpha(x): [0,\max(\alpha,1-\alpha)] \to [0,\max(\alpha,1-\alpha)]$ defined as follows:
\begin{equation*}
M_\alpha(x)=\abs{\frac{1}{x}-\pint{\frac{1}{x}+1-\alpha}}
\end{equation*} 
Notice that for a given $\alpha$, $M_\alpha$ is a factor of $T_\alpha$: in fact $T_{\alpha} \circ h=h \circ M_\alpha$, where $h : x \mapsto \abs{x}$ is the absolute value. Since all the corresponding results for the maps $M_\alpha$ can be derived through this semiconjugacy, in the following paragraphs we will focus on the maps $T_\alpha$.\\
The following proposition extends the results of Moussa, Cassa and Marmi \cite{MCM} to the maps $T_{\alpha}$:
\begin{proposition}
For $\sqrt{2}-1 \leq \alpha \leq \frac{1}{2}$, a representation of the natural extension of $T_\alpha$ is $\overline{T}_{\alpha} : D_{\alpha} \to D_{\alpha}$, where
\begin{align*}
&D_{\alpha}=\left[\alpha-1,\frac{2\alpha-1}{1-\alpha}\right) \times [0,1-g] \cup \\ &
\cup\left[\frac{2\alpha-1}{1-\alpha},\frac{1-2\alpha}{\alpha}\right)\times \left([0,1-g] \cup
\left[\frac{1}{2},g\right]\right) \cup \left[\frac{1-2\alpha}{\alpha},\alpha\right]\times[0,g] \subset \mathbb{R}^2, \notag \\
&\overline{T}_{\alpha} (x,y)=\left(T_{\alpha}(x),\frac{1}{\left[\abs{\frac{1}{x}}+1-\alpha\right]+\sign(x)y}\right)
\end{align*}
Then the invariant density for $\sqrt{2}-1 \leq \alpha \leq \frac{1}{2}$ is
\begin{multline*}
\rho_\alpha(x)=\frac{1}{\log G}\left(\chi_{\left[\alpha-1,\frac{2\alpha-1}{1-\alpha}\right)}(x)\frac{1}{x+G+1}+ \right.\\
\left. +\chi_{\left[\frac{2\alpha-1}{1-\alpha},\frac{1-2\alpha}{\alpha}\right)}\left(\frac{1}{x+G+1}+\frac{1}{x+G}-\frac{1}{x+2}\right) +\chi_{\left[\frac{1-2\alpha}{\alpha},\alpha\right)}(x)\frac{1}{x+G}\right)
\end{multline*}
\end{proposition}      

The proof of the above proposition is quite straightforward and we leave it to the reader.\bigskip
\par It can be shown \cite{DKU} that the Kolmogorov-Sinai entropy with respect to the unique absolutely continuous invariant measure $\mu_{\alpha}$ of the $T_{\alpha}$ is given by Rohlin's formula:
\begin{displaymath}
h(T_{\alpha})=\int_{\alpha-1}^{\alpha} \log{\abs{T_{\alpha}'(x)}} d\mu_{\alpha}(x)
\end{displaymath}
Actually, Rohlin's formula applies also to the $M_\alpha$, and $h(T_{\alpha})=h(M_{\alpha})$.
For $\sqrt{2}-1 \leq \alpha \leq 1$, the entropy can be computed explicitly from the expression of the invariant densities \cite{Nak}, \cite{MCM}:
\begin{align} \label{entropy2}
h(T_{\alpha})=
\left\{\begin{array}{ll} 
\frac{\pi^2}{6 \log(1+\alpha)} & \text{for } g < \alpha \leq 1 \\
\frac{\pi^2}{6\log G} & \text{for } \sqrt{2}-1 \leq \alpha \leq g  
\end{array}
\right.
\end{align}  
In particular, the entropy is constant when $\sqrt{2}-1 \leq \alpha \leq g$ and its derivative has a discontinuity (\emph{phase transition}) in $\alpha=g$.\bigskip
\par The case $\alpha=0$ requires a separate discussion; in fact, due to the presence of an indifferent fixed point, $T_0$ doesn't admit a finite invariant density, although it is invariant with respect to the infinite measure $d\mu_0=\frac{dx}{1+x}$. Therefore the entropy of $T_0$ can only be defined in \emph{Krengel's sense}, that is up to multiplication by a constant (see M. Thaler \cite{Th} for a study of the general one-dimensional case).
Following \cite{Th}, for any subset $A$ of $[0,1]$ with $0< \mu_0(A) <\infty$ we can define
\begin{equation*} 
h(T_0,\mu_0) \doteqdot \mu_0(A) h(({T_0})_A)
\end{equation*} 
where $h(({T_0})_A)$ is the entropy of the first return map of $T_0$ on $A$ with respect to the normalized induced measure $\mu_A=\frac{\mu_0}{\mu_0(A)}$. This quantity is well-defined since the product $h(T_0,\mu_0)$ doesn't depend on the choice of $A$, and it has been computed exactly: $h(T_0,\mu_0)=\frac{\pi^2}{3\log2}$ \cite{Va}. Since this is a finite value, for a sequence $A_k$ of subsets whose Lebesgue measure tends to $1$ we would have $h(({T_0})_{A_{k}})=\frac{\pi^2}{(3\log2)\mu_0(A_k)}\to 0$. In this restricted sense we can say that \textquotedblleft the entropy of $T_0$ is $0$\textquotedblright.\bigskip
\par Expression (\ref{entropy2}) suggests the notion that the dynamical systems $T_\alpha$ are somehow related and have a common origin; actually for $\frac{1}{2} \leq \alpha \leq g$ their natural extensions are all isomorphic. Moreover, a recent result by R. Natsui \cite{Nat} shows that the natural extensions of the Farey maps associated to the $T_\alpha$ are all isomorphic when $\frac{1}{2} \leq \alpha \leq 1$.\\
It is well-known that the maps $T_1$ and $T_0$ descend from the geodesic flow on the unit tangent bundle of the modular surface $PSL(2,\mathbb{Z})\backslash PSL(2,\mathbb{R})$ \cite{Se}, \cite{Ka}. Indeed we can represent this flow as a suspension flow over the natural extension of these maps and deduce in this way the invariant probability measures from the normalized Haar measure on $PSL(2,\mathbb{Z})\backslash PSL(2,\mathbb{R})$. It is natural to conjecture that the same happens for all the maps $T_\alpha,\, \alpha \in [0,1]$. If this were true, one could (at least in principle) apply Abramov's formula to compute the entropies $h(\alpha)$ from the entropy of the geodesic flow.\\ 

We now summarize briefly the contents of the various sections of the paper.\\
In \S\, \ref{continuity_section} we prove that the entropy $h(\alpha)$ of $T_{\alpha}$ is continuous in $\alpha$ when $\alpha \in (0,1]$ and that $h(\alpha) \to 0$ as $\alpha \to 0$, as it had been conjectured by Cassa \cite{Ca}. This result is based on a uniform version of the Lasota-Yorke inequality for the Perron-Frobenius operator of $T_{\alpha}$, following M. Viana's approach \cite{Vi}; in the uniform case, however, a further difficulty arises from the existence of arbitrarily small cylinders containing the endpoints, requiring \emph{ad hoc} estimates.\\
In \S\, \ref{numerical_simulations} we analyse the results of numerical simulations for the entropy obtained through Birkhoff sums, which suggest that the entropy function has a complex self-similar structure.\\ Finally, in \S\, \ref{NE} we compute the natural extension and the invariant densities of the $T_{\alpha}$ for the sequence $\left\{\alpha=\frac{1}{r}\right\}_{r \in \mathbb{N}}$.

\section*{Acknowledgements}
We are grateful to Viviane Baladi and Marcelo Viana for several useful discussions. In particular we wish to thank the latter for suggesting to use the jump transformation in the proof of Proposition \ref{delta}. We are also grateful to Rie Natsui for explaining us her work on $\alpha$-Farey maps.\\
One of us (L.L.) would like to thank Roberto Pinciroli for his careful reading of the manuscript, including many helpful suggestions.\\ 
We wish to acknowledge the financial support of the MURST project \textquotedblleft Sistemi dinamici non lineari e applicazioni fisiche\textquotedblright\, of the Scuola Normale Superiore and of the Centro di Ricerca Ma\-te\-ma\-ti\-ca \textquotedblleft Ennio De Giorgi\textquotedblright.

\section{Continuity of the entropy} \label{continuity_section}
The main goal of the present section is the following
\begin{theorem}\label{main}
The function $\alpha \to h(\alpha)$ is continuous in $(0,1]$, and 
\begin{displaymath}
\lim_{\alpha\to 0^+} h(\alpha) =0
\end{displaymath}
\end{theorem}

Since in the case $\alpha \geq \sqrt{2}-1$ the entropy has been computed exactly by Nakada \cite{Nak} and Marmi, Moussa, Cassa \cite{MCM}, we can restrict our study to the case  $0<\alpha \leq \sqrt{2}-1$.

It is well-known \cite{BDV} that for all $\alpha \in (0,1]$ the maps $T_\alpha$ admit a unique absolutely continuous invariant probability measure $\mu_\alpha$, whose density $\rho_\alpha$ is of bounded variation (and therefore bounded). In addition, a result of R. Zweim\"uller entails that $\rho_\alpha$ is bounded from below (see \cite{Zw1}, Lemma 7): 
\begin{equation} \label{C}
\forall \alpha \in (0,\sqrt{2}-1],\; \exists C>0 \text{ s.t }\forall x \in I_{\alpha}, \rho_{\alpha}(x)\geq C
\end{equation}
The uniqueness of the a.c.i.m. is a consequence of the ergodicity of the system: 

\begin{lemma}[Exactness] \label{exactness}
For all $\alpha \in [0,1]$, the dynamical system $(T_{\alpha},\mu_{\alpha})$ is exact (and therefore ergodic).
\end{lemma}

The proof of this Lemma for $\alpha \in \left[\frac{1}{2},1\right]$ was given by H. Nakada \cite{Nak} and can be adapted to our case with slight changes (see Appendix \ref{exact}). The fact that $T_0$ is exact follows from a result of M. Thaler \cite{Th}.\bigskip

\par To prove continuity we adopt the following approach: by means of a uniform Lasota-Yorke-type inequality for the Perron-Frobenius operator, we prove that the variations of the invariant densities are equibounded as $\alpha$ varies in some neighborhood of any fixed $\bar{\alpha}>0$ (see Proposition \ref{BV} below). Our argument follows quite closely \cite{Vi}, except that we have to deal with a further difficulty arising from the fact that the cylinders containing the endpoints $\alpha$ and $\alpha-1$ can be arbitrarily small. After translating the maps so that their interval of definition does not depend on $\alpha$ around $\bar{\alpha}$, we prove the $L^1$-continuity of the invariant densities $\rho_{\alpha}$ using Helly's Theorem (Lemma \ref{continuity}). Then the continuity of the entropy follows from Rohlin's formula.  

\subsection{Notations}

\subsubsection{Cylinders of rank $1$}
Let $0< \alpha \leq \sqrt{2}-1$. The map $T_\alpha$ is piecewise analytic on the countable partition $\mathcal{P}=\{I^{+}_j\}_{j\geq j_{\min}} \cup \{I^{-}_{j}\}_{j \geq 2}$, 
where $j_{\min}=\left[\abs{\frac{1}{\alpha}}+1-\alpha\right]$, \quad and the elements of $\mathcal{P}$ are called \emph{cylinders of rank $1$}:  
\begin{align*}
& I^{+}_{j}\doteqdot \left(\frac{1}{j+\alpha},\frac{1}{j-1+\alpha}\right], \; j \in [j_{\min}+1,\infty),\quad	I^{+}_{j_{\min}}\doteqdot\left(\frac{1}{j_{\min}+\alpha},\alpha\right],\\
& I^{-}_{j}\doteqdot\left[-\frac{1}{j-1+\alpha},-\frac{1}{j+\alpha}\right), \; j \in [3,\infty), \quad
I^{-}_{2}\doteqdot\left[\alpha-1,-\frac{1}{2+\alpha}\right)
\end{align*}
$T_{\alpha}$ is monotone on each cylinder and we have
\begin{align*}
&\left\{ \begin{array}{ll}
T_{\alpha}(x)=
\frac{1}{x}-j,  & x \in I^{+}_{j},\; j\in \mathbb{N} \cap [j_{\min},\infty) \\
T_{\alpha}(x)= -\frac{1}{x}-j, &x \in I^{-}_{j},\; j \in \mathbb{N} \cap [2,\infty)
\end{array} \right.  \\
\end{align*}
Thus for $x \in I^{\pm}_{j}$,
\begin{align} \label{l}
&\frac{1}{\abs{T_{\alpha}'(x)}}\leq \lambda_j \leq \lambda <1\notag , \intertext{ where } 
&\lambda= (1-\bar{\alpha}+\varepsilon)^2, \quad \lambda_j=\frac{1}{(j-1+\bar{\alpha}-\varepsilon)^2},\; j>2,\quad \lambda_2=\lambda
\end{align}
depend only on $\bar{\alpha}$ and $\varepsilon$. Moreover, we have that $\Var_{I^{\pm}_{j}} \abs{\frac{1}{T_{\alpha}'(x)}} \leq \lambda_j \; \forall \alpha \in [\bar{\alpha} - \varepsilon, \bar{\alpha} + \varepsilon]$.
\subsubsection{Cylinders of rank $n$; full cylinders}
Let $\mathcal{P}^{(n)}= \bigvee_{i=0}^{n-1} T_{\alpha}^{-i}(\mathcal{P})$ be the induced partition in monotonicity intervals of $T_{\alpha}^n$. Each cylinder $I_{\eta}^{(n)} \in \mathcal{P}^{(n)}$ is uniquely determined by the sequence 
\begin{equation*}
((j_0(\eta),\varepsilon_0(\eta)),\ldots,(j_{n-1}(\eta),\varepsilon_{n-1}(\eta))
\end{equation*}
such that for all $x \in I_{\eta}^{(n)}, \; T_{\alpha}^i(x) \in I_{j_i(\eta)}^{\varepsilon_i(\eta)}$. On each cylinder $T_{\alpha}^n$ is a M\"obius map $T_{\alpha}^n(x)=\frac{ax+b}{cx+d}$, where $\left(\begin{array}{cc} a & b \\ c & d\end{array}\right) \in GL(2,\mathbb{Z})$. We will say that a cylinder $I_{\eta}^{(n)} \in \mathcal{P}$ is \emph{full} if $T_{\alpha}^n(I_{\eta}^{(n)})=I_{\alpha}$.

\subsubsection{Perron-Frobenius operator}
Let $V_\eta:T_{\alpha}^n(I_{\eta}^{(n)}) \to I_{\eta}^{(n)}$ be the inverse branches of $T_{\alpha}^n$, and $P_{T_{\alpha}}$ the Perron-Frobenius operator associated with $T_{\alpha}$. Then for every $\varphi \in L^1(I_{\alpha})$, 
\begin{align}
(P_{T_{\alpha}}^{n} \varphi)(x)= \sum_{I_{\eta}^{(n)} \in \mathcal{P}^n} \frac{\varphi\left(V_{\eta}(x)\right)}{\abs{(T_{\alpha}^n)'\left(V_{\eta}(x)\right)}} \chi_{T_{\alpha}^n(I_{\eta}^{(n)})}(x)  \label{P_n}
\end{align}

On $I_{\eta}^{(n)}$ we have the following bound:  
\begin{equation*}
\sup_{I_{\eta}^{(n)}}\frac{1}{\abs{\left(T_{\alpha}^n\right)'(x)}}=\sup_{I_{\eta}^{(n)}} \frac{1}{\abs{T_{\alpha}'\left(T_{\alpha}^{n-1}(x)\right)\cdots T_{\alpha}'(x)}} \leq \lambda_{\eta}^{(n)} \leq \lambda^n,
\end{equation*}
where $\lambda_{\eta}^{(n)} \doteqdot \lambda_{{j_0}(\eta)}\cdots \lambda_{{j_{n-1}}(\eta)}$.  
Recall that for $f_1,\ldots,f_n \in BV$,
\begin{equation} \tag{i} \label{var_prod}
\Var(f_1\cdots f_n) \leq \sum_{k=1}^{n} \Var(f_k) \prod_{i \neq k} \sup \abs{f_i}
\end{equation} 
and consequently
\begin{align*} 
&\Var_{I_{\eta}^{(n)}} \frac{1}{\abs{\left(T_{\alpha}^n\right)'(x)}}=\Var_{I_{\eta}^{(n)}} \frac{1}{\abs{T_{\alpha}'\left(T_{\alpha}^{n-1}(x)\right)\cdots T_{\alpha}'(T_{\alpha}(x))\cdot T_{\alpha}'(x)}}  \leq  n \lambda_{\eta}^{(n)}
\end{align*}

\subsubsection{}
Finally, we state the following \emph{bounded distortion property}, that we are going to use several times in the sequel:
\begin{proposition}[Bounded distortion] \label{bound_dist}
$\forall \alpha >0$, $\exists C_1$ such that $\forall n \geq 1, \forall I_{\eta}^{(n)} \in \mathcal{P}^{(n)}, \forall x,y \in I_{\eta}^{(n)}$,
\begin{equation*}
\abs{\frac{(T^n_{\alpha})'(y)}{(T^n_{\alpha})'(x)}} \leq C_1
\end{equation*}
Moreover, for all measurable set $B \subseteq I_\alpha$, for all full cylinders $I_{\eta}^{(n)} \in \mathcal{P}^{(n)}$,
\begin{equation*}
m(V_\eta(B)) \geq \frac{m(B)m(I_{\eta}^{(n)})}{C_1},
\end{equation*}
where $m$ denotes the Lebesgue measure.
\end{proposition}  
The proof of this statement follows a standard argument and can be found in Appendix \ref{bounded_distortion}.

\subsection{Uniformly bounded variation of the invariant densities} \label{UBV}
Let $\bar{\alpha} \in (0,\sqrt{2}-1] $ and $\varepsilon<\bar{\alpha}$ be fixed, and choose $\alpha \in [\bar{\alpha}-\varepsilon,\bar{\alpha} + \varepsilon]$. 
     
\begin{proposition} \label{BV}
$\forall \bar{\alpha} \in (0,\sqrt{2}-1]$, $\rho_{\bar{\alpha}}$ is of bounded variation, and $\exists \varepsilon$, $\exists K>0$ such that for all $\alpha \in [\bar{\alpha}-\varepsilon,\bar{\alpha}+\varepsilon]$, $\Var(\rho_\alpha)<K$.    
\end{proposition}

The main result we need in order to prove Proposition \ref{BV} is the following

\begin{lemma}[Uniform version of the Lasota-Yorke inequality] \label{Lemma}
Let $\bar{\alpha}$ be fixed. Then there exist $\lambda_0<1$, $C,K_0>0$ such that $\forall n$, $ \forall \varphi \in BV(I)$, $\forall \alpha \in [\bar{\alpha}-\varepsilon,\bar{\alpha}+\varepsilon]$, 
\begin{equation*}
\Var_{I_\alpha} \left(P_{T_{\alpha}}^{n} \varphi \right) \leq C (\lambda_0)^n \Var \varphi + K_0 \int\limits_{I_\alpha} \abs{\varphi} dx
\end{equation*}
\end{lemma}

Assuming Lemma \ref{Lemma} the Proposition then follows easily. Indeed it is enough to recall that the Cesaro sums  
\begin{equation*}
\rho_n= \frac{1}{n} \sum_{j=0}^{n-1} P_{T_{\alpha}}^j 1
\end{equation*}
of the sequence $\{P_{T_{\alpha}}^{j} 1\}_{j \in \mathbb{N}}$ converge almost everywhere to the invariant density ${\rho}_{\alpha}$ of $T_\alpha$. Both the variations and the $L^{\infty}$ norms of the $\{\rho_n\}$ are uniformly bounded: 
\begin{align*}
&\Var \rho_n \leq \frac{1}{n} \sum_{j=0}^{n-1} \Var \left(P_{T_{\alpha}}^j 1 \right) \leq \frac{1}{n} \sum_{j=0}^{n-1} K_0 m(I_\alpha) =K_0 \quad \forall n \\
&\int\limits_{I_\alpha} \rho_n dx = \frac{1}{n} \sum_{j=0}^{n-1} \int\limits_{I_\alpha} P_{T_{\alpha}}^j 1 dx = m(I_{\bar{\alpha}}) = 1 \quad \forall n \Rightarrow \\
&\sup_{I_\alpha} \abs{\rho_n} \leq \Var_{I_\alpha} \rho_n + \frac{1}{m(I_\alpha)} \leq K_0 +1 \quad \forall n, 
\end{align*}
where $K_0$ is the constant we found in Lemma \ref{Lemma}. Then we also have $\Var {\rho}_{\alpha} \leq K_0,\; \sup \abs{\rho_{\alpha}} \leq K_0+1$, which concludes the proof of Proposition \ref{BV}.  

\begin{proof}[\textbf{Proof of Lemma \ref{Lemma}}] 
We have
\begin{align}
& \Var \left(P_{T_{\alpha}}^{n} \varphi \right) \leq \sum_\eta \left( \Var_{T_{\alpha}^n(I_{\eta}^{(n)})} \frac{\varphi\left(V_{\eta}(x)\right)}{\abs{(T_{\alpha}^n)'\left(V_{\eta}(x)\right)}} + 2 \sup_{T_{\alpha}^n(I_{\eta}^{(n)})} \abs{\frac{\varphi\left(V_{\eta}(x)\right)}{(T_{\alpha}^n)'\left(V_\eta(x)\right)}}  \right) = \notag \\
& = \sum_\eta \left( \Var_{I_{\eta}^{(n)}} \frac{\varphi(y)}{\abs{(T_{\alpha}^n)'(y)}} + 2 \sup_{I_{\eta}^{(n)}} \abs{\frac{\varphi\left(y\right)}{(T_{\alpha}^n)'\left(y\right)}}  \right) \label{P_n2}
\end{align}
For the last equality, observe that since $V_\eta:T_{\alpha}^n(I_{\eta}^{(n)}) \to I_{\eta}^{(n)}$ is a homeomorphism, $\Var_{T_{\alpha}^n(I_{\eta}^{(n)})} \left(\frac{\varphi}{\abs{(T_{\alpha}^n)'}} \circ V_\eta\right)= \Var_{I_{\eta}^{(n)}} \frac{\varphi}{\abs{(T_{\alpha}^n)'}}$. 
The first term in expression (\ref{P_n2}) can be estimated using (\ref{var_prod}):
\begin{align*}
&\sum_\eta \Var_{I_{\eta}^{(n)}} \frac{\varphi(y)}{\abs{(T_{\alpha}^n)'(y)}} \leq \sum_\eta \left( \Var_{I_{\eta}^{(n)}} \varphi \sup_{I_{\eta}^{(n)}}\frac{1}{\abs{(T_{\alpha}^n)'(y)}} + \Var_{I_{\eta}^{(n)}} \frac{1}{\abs{(T_{\alpha}^n)'(y)}}\sup_{I_{\eta}^{(n)}} \abs{\varphi} \right) \leq \\
& \leq \sum_\eta \left( \lambda_{\eta}^{(n)} \Var_{I_{\eta}^{(n)}} \varphi + n \lambda_{\eta}^{(n)} \sup_{I_{\eta}^{(n)}} \abs{\varphi} \right) 
\end{align*}
For the second term, we have $2 \sum_\eta \sup_{I_{\eta}^{(n)}} \abs{\frac{\varphi\left(y\right)}{(T_{\alpha}^n)'\left(y\right)}} \leq 2 \sum_\eta \lambda_{\eta}^{(n)} \sup_{I_{\eta}^{(n)}} \abs{\varphi(y)}$. In conclusion, from equation (\ref{P_n2}) we get 
\begin{equation} \label{var_P}
\Var \left(P_{T_{\alpha}}^n \varphi \right)(x) \leq \lambda^n \Var_{I_\alpha} \varphi + \sum_{\eta} (n +2) \lambda_{\eta}^{(n)} \sup_{I_{\eta}^{(n)}}\abs{\varphi}  
\end{equation} 
We want to give an estimate of the sum in equation (\ref{var_P}); recall that for $\varphi \in BV$,
\begin{equation} \tag{ii} \label{sup_int}
\sup_{I_{\eta}^{(n)}} \abs{\varphi} \leq \Var_{I_{\eta}^{(n)}} \varphi + \frac{1}{m(I_{\eta}^{(n)})}\int\limits_{I_{\eta}^{(n)}} \abs{\varphi} dx
\end{equation}
However, equation (\ref{sup_int}) doesn't provide a global bound independent from $\eta$ for two reasons. In the first place, the lengths of the intervals $I_{\eta}^{(n)}$ are not bounded from below when the indices $j_{i}(\eta)$ grow to infinity. Furthermore, a difficulty that arises only in the case of uniform continuity and that was not dealt with in reference \cite{Vi} is that the measures of the cylinders of rank $n$ containing the endpoints $\bar{\alpha}$ and $\bar{\alpha}-1$ are \emph{not} uniformly bounded from below in $\alpha$, and require a careful handling.\\ 
To overcome the first difficulty, following \cite{Vi}, we split the sum into two parts: for $n$ fixed, let $k$ be such that
\begin{equation} \label{k}
\sum_{j>k} \lambda_j \leq \frac{\lambda^n}{2^{2n-1}}
\end{equation}
Since $\lambda$ doesn't depend on $\alpha$, neither does $k$. 
Define the set of \textquotedblleft intervals with bounded itineraries'' 
\begin{equation} \label{G(n)}
G(n)=\{ I_{\eta}^{(n)} \in \mathcal{P}^n \;|\; \max(j_0(\eta),\ldots,j_{n-1}(\eta)) \leq k\}
\end{equation}
To get rid of the measures of the cylinders containing the endpoints, we combine them with full cylinders; the measures of the latter can be estimated using Lagrange's Theorem, since the derivatives are bounded under the hypothesis of bounded itineraries. When combining intervals, we have to consider the sum of the corresponding $\lambda_{v}^{(n)}$ and make sure that it is smaller than $1$. This requires additional care when $I_{\eta}^{(n)} \ni \alpha-1$.  
\begin{remark} \label{nu} 
Let $r=r(\alpha)$ be such that 
\begin{equation} \label{v_r}
v_{r+1} \leq \alpha < v_r, \quad \text{ where } v_{r} =-\frac{1}{2}+\frac{1}{2}\sqrt{1+\frac{4}{r}}
\end{equation}
(clearly $r$ is bounded by $r(\bar{\alpha})+1$ in a small neighborhood of $\bar{\alpha}$). Then  $T_{\alpha}^i(\alpha-1)=\frac{(i+1)\alpha-1}{1-i\alpha}\in I_{2}^{-}$ for $i=0,\ldots,r-1$ and $T_{\alpha}^r(\alpha-1) \notin I_{2}^{-}$. Thus any cylinder with more than $r$ consecutive digits \textquotedblleft$(2,-)$" is empty, and the cylinder $((2,-),\ldots,(2,-))$ of rank $r$ may be arbitrarily small when $\alpha$ varies. The cylinder $(j_{\min},+)$ can be arbitrarily small too.
\end{remark} 
Consider the function $\sigma: G(n) \to G(n)$ which maps every nonempty cylinder $I_{\eta}^{(n)}$ in $I_{\xi}^{(n)}$ in the following way:
\begin{enumerate}
\renewcommand{\theenumi}{\alph{enumi}}
\item If $(j_i(\eta),\varepsilon_i(\eta))=(j_{\min},+)$ for some $i$, then $(j_i(\xi),\epsilon_i(\xi))=(j_{\min}+1,+)$;
\item If $\exists i $ such that 
\begin{displaymath}
((j_i(\eta),\varepsilon_i(\eta)),\ldots,(j_{i+r}(\eta),\varepsilon_{i+r}(\eta)))=\\((2,-),(2,-),\ldots,(2,-)),
\end{displaymath}
then $((j_i(\xi),\varepsilon_i(\xi)),\ldots,(j_{i+r}(\xi),\varepsilon_{i+r}(\xi)))=((2,-),\ldots,(2,-),(3,-))$;
\item Otherwise, $(j_i(\xi),\varepsilon_i(\xi))=(j_i(\eta),\varepsilon_i(\eta))$.
\end{enumerate}
We want to show that there exists $\delta_n>0$, depending only on $\bar{\alpha}$, such that for all $\xi \in \sigma(G(n))$, $m(\xi) \geq \delta_n$. For this purpose, we group together the sequences of consecutive digits $(2,-)$, and obtain a new alphabet $\mathcal{A}=\mathcal{A}_1 \cup \mathcal{A}_2$, where
\begin{align*}
&\mathcal{A}_1=\{(3,-),\ldots,(k,-)\}\cup \{(j_{\min}+1,+),\ldots,(k,+)\}\\
&\mathcal{A}_2=\{(2,-),((2,-),(2,-)),\ldots,(\underbrace{(2,-),\ldots,(2,-))}_{r-1}\}
\end{align*}
Then each $\xi \in \sigma(G(n))$ can be seen as a sequence in $\mathcal{A'}^s=\{(a_1,\ldots,a_s)\in \mathcal{A}^s \,|\, a_i \in \mathcal{A}_2 \Rightarrow a_{i+1} \in \mathcal{A}_1\}$ for some $n \geq s \geq \frac{n}{r}$. Let $\widetilde{T_{\alpha}}$ be the first return map on $\mathcal{A}_1$ restricted to $\sigma(G(n))$:
\begin{equation*}
\begin{array}{ll}
\widetilde{T}(x)=T_{\alpha}(x) & \text{for } x \in (a) \in \mathcal{A}_1; \\
\widetilde{T}(x)=T_{\alpha}^i(x) & \text{if } \exists i: x \in \underbrace{((2,-),\ldots,(2,-))}_{i}, x \notin \underbrace{((2,-),\ldots,(2,-))}_{i+1}
\end{array}
\end{equation*}  
Let $\widetilde{V}_{a}$ be the inverse branch of $\widetilde{T}$ relative to the cylinder $(a)$. Observe that 
\begin{equation} \label{star}
\forall (a_1,\ldots,a_s) \in \mathcal{A}'^{s}, \widetilde{T}^s(a_1,\ldots,a_s)= \widetilde{T}(a_s)
\end{equation}
This can be proved by induction on $s$: when $s=1$ it is trivial; supposing that the property (\ref{star}) holds for all sequences of length $s$, we have 
\begin{multline*}
\widetilde{T}^{s+1}(a_1,\ldots,a_{s+1})=\widetilde{T}^{s+1}\left((a_1,\ldots,a_s)\cap \widetilde{V}_{a_1}\cdots \widetilde{V}_{a_s}(a_{s+1})\right)=\\
=\widetilde{T}(\widetilde{T}^s(a_1,\ldots,a_s)\cap (a_{s+1})) 
\end{multline*} 
since $\widetilde{T}^s$ is injective on $(a_1,\ldots,a_s)$; this is equal to $\widetilde{T}(\widetilde{T}(a_s) \cap (a_{s+1}))$ by inductive hypothesis.  
\begin{itemize}
\item If $a_{s+1} \in \mathcal{A}_2$, we have $a_s \in \mathcal{A}_1$ and $\widetilde{T}(a_s)=I$: then $\widetilde{T}^{s+1}(a_1,\ldots,a_{s+1})=\widetilde{T}(a_{s+1})$. 
\item If $a_{s+1}\in \mathcal{A}_1$, $\widetilde{T}(a_s) \supseteq (a_{s+1})$. In fact for all $i=0,\ldots,r-1$, 
\begin{multline}
T_{\alpha}^i(\underbrace{(2,-),\ldots,(2,-)}_{i})=T_{\alpha}^i\left(\left[\alpha-1,V_{(2,-)}^{i}(\alpha)\right)\right)=\\
=[T_{\alpha}^i(\alpha-1),\alpha)\supseteq \left[-\frac{1}{2+\alpha},\alpha\right)\supseteq \bigcup\limits_{a \in \mathcal{A}_1} (a) \supseteq \left[-\frac{1}{3},0\right] \label{surj}
\end{multline}
\end{itemize}
Equation (\ref{star}) provides a lower bound on the measures of the intervals in $\sigma(G(n))$:
\begin{align*}
&\frac{1}{3} \leq m(\widetilde{T}(a_s))=m(\widetilde{T}^s(a_1,\ldots,a_s))\leq m(a_1,\ldots,a_s) \sup \abs{(\widetilde{T}^s)'}, \;\text{ and }\\
&M(\bar{\alpha})\doteqdot \left(\max\left((k+\alpha+\varepsilon)^{2},(2+\alpha+\varepsilon)^{2r(\bar{\alpha})}\right)\right)\geq \sup \abs{\widetilde{T}'}
\end{align*}
in a neighborhood of $\bar{\alpha}$. Thus for all $I_{\xi}^{(n)} \in \sigma(G(n))$,
\begin{equation*}
\delta_n \doteqdot \frac{1}{3M(\bar{\alpha})^n} \leq m(I_{\xi}^{(n)})
\end{equation*}
Returning to the sum in equation (\ref{var_P}), and defining $\overline{I}_{\xi}^{(n)}=\bigcup\{I_{\eta}^{(n)} \,|\, \sigma(I_{\eta}^{(n)})=I_{\xi}^{(n)}\}$, we find:
\begin{align*}
&\sum_{I_{\eta}^{(n)} \in G(n)} \left( \lambda_{\eta}^{(n)} \sup_{I_{\eta}^{(n)}}\abs{\varphi} \right)  
\leq \sum_{I_{\xi}^{(n)} \in \sigma(G(n))} \sup_{\overline{I}_{\xi}^{(n)}} \abs{\varphi} \left(\sum_{\sigma(I_{\eta}^{(n)})=I_{\xi}^{(n)}}\lambda_{\eta}^{(n)}\right)  
\end{align*}
We want to estimate $\lambda'=\sup\limits_{\sigma(G(n))} \sum\limits_{\sigma(I_{\eta}^{(n)})=I_{\xi}^{(n)}}\lambda_{\eta}^{(n)}$:
each sum can be computed distributively as a product of at most $n$ factors $\lambda'_i$, each of which corresponds to one of the cases a),\,b),\,c) that we have listed in the definition of $\sigma$:
\begin{itemize}
\item In the case a), we have $\lambda'_i=\lambda_{j_{\min}}+\lambda_{j_{\min}+1} \leq  2(\bar{\alpha}+\varepsilon)^2<\frac{1}{2}$ (remark that $j_{\min}\geq 3$ when $\alpha \leq \sqrt{2}-1$).
\item In the case b), $\lambda'_i=\lambda_{2}^{r}+\lambda_{2}^{r-1}\lambda_{3}=(1-\alpha)^{2(r-1)}\left((1-\alpha)^2+\frac{1}{(2+\alpha)^2}\right)<0.9$. In fact, when $\alpha>\frac{1}{5}$ we have $(1-\alpha)^2+\frac{1}{(2+\alpha)^2}<\frac{9}{10}$; otherwise,  $(1-\alpha)^2+\frac{1}{(2+\alpha)^2}<\frac{5}{4},$ and for $\alpha \geq \eta_{r+1}$, we have $r-1 \geq\frac{1}{\alpha^2+\alpha}-2$, and
\begin{align*}
&(1-\alpha)^{2(r-1)}=\left(\frac{1-\alpha^2}{1+\alpha}\right)^{2(r-1)} \leq \frac{1}{(1+\alpha)^{2(r-1)}}\leq \frac{1}{1+2\alpha(r-1)} \leq \\
&\leq \frac{1+\alpha}{3-3\alpha-4\alpha^2} < \frac{3}{5}
\end{align*}
\item In the case c), $\lambda'_i=\lambda_{j_i}$. 
\end{itemize} 
(The constants in the previous discussion are far from optimal, but they are sufficient for our purposes.)\\ 
Then $\lambda' \leq \max\left(\lambda^n,\left(\frac{9}{10}\right)^{\frac{n}{r(\bar{\alpha})+1}}\right)= \tilde{\lambda}^n<1$. Note that $\tilde{\lambda}$ only depends on $\bar{\alpha}$ and not on $\alpha$.\\
We can finally complete our estimate for the sum over $I_{\eta}^{(n)} \in G(n)$:
\begin{align} \label{tilde}
&\lambda' \sum_{I_{\xi}^{(n)} \in \sigma(G(n))} \sup_{\overline{I}_{\xi}^{(n)}} \abs{\varphi} 
\leq \tilde{\lambda}^n \sum_{I_{\xi}^{(n)} \in \sigma(G(n))} \left(\Var_{\overline{I}_{\xi}^{(n)}} \abs{\varphi}+ \frac{1}{m(\overline{I}_{\xi}^{(n)})}\int_{\overline{I}_{\xi}^{(n)}} \varphi \right) \leq \notag \\
&\leq \tilde{\lambda}^n \left(\Var \varphi + \sum_{I_{\xi}^{(n)} \in \sigma(G(n))} \frac{1}{m(I_{\xi}^{(n)})} \int_{\overline{I}_{\xi}^{(n)}} \varphi\right) \leq \tilde{\lambda}^n \Var \varphi + \frac{\tilde{\lambda}^n}{\delta_n} \norm{\varphi}_1  
\end{align} 
On the other hand, for the sum over $I_{\eta}^{(n)} \notin G(n)$ we have the following estimate:
\begin{align}
&\sum_{I_{\eta}^{(n)} \notin G(n)} \left((n+2) \lambda_{\eta}^{(n)} \sup_{I_{\eta}^{(n)}}\abs{\varphi} \right) \leq \notag \\
& \leq \sup_{I_\alpha} \abs{\varphi} \sum_{j>k} \sum_{l=0}^{n-1} \sum_{j_{l}(\eta)=\max \{j_0(\eta),\ldots,j_{n-1}(\eta)\}=j} (n+2) \lambda_{j_0(\eta)} \cdots \lambda_{j_{n-1}(\eta)} \label{9} 
\end{align}
where in the third sum of expression (\ref{9}) we take $l$ to be the smallest integer that realizes the maximum, to avoid counting the same sequences twice. Observe that when we take the sum over $j_0(\eta),\ldots,j_{n-1}(\eta)$, since we are not taking into account the signs $\varepsilon_i(\eta)$, we are actually counting at most $2^n$ distinct sequences.
\begin{align*}
&\sum_{\substack{(j_0(\eta),\ldots,j_{n-1}(\eta)) \\ j_l(\eta)=j}} \lambda_{j_0(\eta)} \cdots \lambda_{j_{n-1}(\eta)} \leq \\
&\leq \lambda_j \left(2 \sum_{(j_0(\eta),\ldots,j_{l-1}(\eta),j_{l+1}(\eta),\ldots,j_{n-1}(\eta))} \lambda_{j_{0}(\eta)} \cdots \lambda_{j_{n-1}(\eta)}\right) \leq \\
&\leq \lambda_j \left( 2\prod_{\substack{i=0 \\ i \neq l}}^ {n-1} \sum_{j_i=2}^j 4\lambda_{j_i} \right) \leq \lambda_j {2}^{2n-1}
\end{align*}
since $\sum_2^{\infty}{\lambda_j}\leq \frac{\pi^2}{6}\leq 2$. Therefore
\begin{align*}
&\sum_{I_{\eta}^{(n)} \notin G(n)} \left((n +2) \lambda_{\eta}^{(n)} \sup_{I_{\eta}^{(n)}}\abs{\varphi} \right) \leq \sup_{I_\alpha} \abs{\varphi} \sum_{j>k} \sum_{l=0}^{n-1} (n +2) \lambda_j 2^{2n-1}\leq \\
& \leq \sup_{I_\alpha} \abs{\varphi}  n(n +2) 2^{2n-1} \sum_{j>k} \lambda_j  \leq \sup_{I_\alpha} \abs{\varphi} n(n + 2) \lambda^n
\end{align*}
where in the last inequality we have used the hypothesis (\ref{k}) on $k$.\\
In conclusion, $\Var\limits_{I_\alpha}(P_{T_{\alpha}}^n \varphi)$ is bounded by
\begin{align*}
\tilde{\lambda}^n \left( (n^2+3n+3) \Var_{I_\alpha}\varphi+(n +2)\left(\frac{1}{\delta_n}+n\right) \norm{\varphi}_1  \right)
\end{align*}
and we recall that we have chosen $\delta_n$ and $\tilde{\lambda}$ so that they do not depend on $\alpha$. \\
Choose any $\bar{\lambda} \in (\tilde{\lambda},1)$, and let $\overline{K}>0, N \in \mathbb{N}$ be such that 
\begin{equation*}
\forall n \geq 1, \;(n^2+3n+3) {\tilde{\lambda}}^n \leq \overline{K} {\bar{\lambda}}^n \text{ and }  \forall n \geq N, \;  \overline{K} {\bar{\lambda}}^n \leq \frac{1}{2} 
\end{equation*}
Let $L(n)= (n+2)\left(\frac{1}{\delta_n}+n\right)\bar{\lambda}^n$, $\hat{K}=\max\limits_{1 \leq n \leq N} L(n)$. For any $n$, we can perform the Euclidean division $n=qN+r$ for some $q \geq 0$ and $0 \leq r < N$.  
Then 
\begin{align}
&\Var_{I_{\alpha}} \left(P_{T_{\alpha}}^N \varphi\right) \leq 
\overline{K} \bar{\lambda}^N \Var_{I_\alpha} \varphi + \hat{K}\norm{\varphi}_1 \label{1a}
\end{align}
More generally, we can show by induction on $q$ that 
\begin{equation} \label{induction}
\Var_{I_\alpha} \left(P_{T_{\alpha}}^{qN}\varphi \right) \leq (\overline{K} \lambda^N)^q \Var_{I_\alpha} \varphi + C(q) \hat{K} \norm{\varphi}_1
\end{equation}
where $C(q)=1+ \frac{1}{2}+\cdots+ \frac{1}{2^{q-1}} <2$ for all $q$. In fact if (\ref{induction}) is true for some $q$, recalling that the Perron-Frobenius operator $P_{T_{\alpha}}$ preserves the $L^1$ norm, we get
\begin{align*}
&\Var_{I_\alpha} \left(P_{T_{\alpha}}^{(q+1)N} \varphi \right) \leq (\overline{K} \bar{\lambda}^N)^q \Var_{I_{\alpha}} \left(P_{T_{\alpha}}^{N}\varphi\right) + C(q) \hat{K}  \norm{P_{T_{\alpha}}^{N} \varphi}_1 \leq \\
& \leq (\overline{K} \bar{\lambda}^N)^{q+1} \Var_{I_\alpha} \varphi + \left(C(q)+ \frac{1}{2^q}\right) \hat{K} \norm{\varphi}_1 dx \leq \\
& \leq ( \overline{K} \bar{\lambda}^N)^{q+1} \Var_{I_{\alpha}} \varphi + C(q+1)\hat{K} \norm{\varphi}_1 dx 
\end{align*}
For $0 \leq r<N$, $\Var \left (P_{T_{\alpha}}^{r} \varphi \right) \leq \overline{K} \bar{\lambda}^r  \Var \varphi  + \hat{K} \norm{\varphi}_1$. 
In general, for $n=qN+r$, we obtain 
\begin{align*}
&\Var_{I_\alpha} \left (P_{T_{\alpha}}^{n} \varphi \right) \leq (\overline{K} \bar{\lambda}^N)^q \Var_{I_\alpha} \left(P_{T_{\alpha}}^r \varphi \right) + C(q) \hat{K} \norm{\varphi}_1 \leq  \\
& \leq (\overline{K} \bar{\lambda}^N)^q \overline{K} \bar{\lambda}^r \Var_{I_\alpha} \varphi + \hat{K}\left((\overline{K} \bar{\lambda}^N)^q + C(q) \right) \norm{\varphi}_1 \leq \frac{\overline{K}}{2^{q}} \bar{\lambda}^r \Var_{I_\alpha} \varphi + 3\hat{K} \norm{\varphi}_1  
\end{align*}
Now take $\lambda_0\geq \max \left(\frac{1}{2^{\frac{1}{N}}}, \bar{\lambda} \right)$, so that $\frac{\bar{\lambda}^r}{2^q} \leq (\lambda_0)^r (\lambda_0)^{Nq} = (\lambda_0)^n$. 
This concludes the proof of Lemma \ref{Lemma}. 
\end{proof}

\subsection{$L^1$ continuity of the densities $\rho_\alpha$ and continuity of the entropy}

Let $\bar{\alpha} \in (0, \sqrt{2}-1]$ be fixed. To study the $L^1$-continuity property of the densities $\rho_\alpha$ (and the continuity of the entropy $h(\alpha)$) it is convenient to work with measures supported on the same interval. Thus we rescale the maps $T_\alpha$ with $\alpha$ in a neighborhood of $\bar{\alpha}$ to the interval $[\bar{\alpha}-1,\bar{\alpha}]$ by applying the translation $\tau_{\bar{\alpha}-\alpha}$. Let $A_{\alpha,\bar{\alpha}}=\tau_{\bar{\alpha}-\alpha}\circ T_{\alpha} \circ \tau_{\bar{\alpha}-\alpha}^{-1}$ be the new maps: 
\begin{equation*}
A_{\alpha,\bar{\alpha}}(x)= \abs{\frac{1}{x-\bar{\alpha}+\alpha}}-\left[\abs{\frac{1}{x-\bar{\alpha}+\alpha}}+1-\alpha\right]+\bar{\alpha}-\alpha
\end{equation*}
Let $J^{\pm}_{j}=I^{\pm}_{j}+\bar{\alpha}-\alpha$ be the translated versions of the intervals of the original partition, and $\tilde{\rho}_\alpha(x)=\rho \circ \tau_{\bar{\alpha}-\alpha}^{-1}(x)=\rho(x-\bar{\alpha}+\alpha)$ the invariant densities for $A_{\alpha,\bar{\alpha}}$. Clearly the bounds for the sup and the variation of $\rho_{\alpha}$ are still valid for $\tilde{\rho}_\alpha$.

\begin{lemma}\label{continuity}
Let $\bar{\alpha} \in (0,\sqrt{2}-1]$ be fixed, and let $\varepsilon$ be given by Proposition \ref{BV}. Then if $\{\alpha_n\} \subset [\bar{\alpha}-\varepsilon,\bar{\alpha} + \varepsilon]$ is a monotone sequence converging to $\bar{\alpha}$, we have $\tilde{\rho}_{\alpha_n} \xrightarrow{\;L^1\;} \tilde{\rho}_{\bar{\alpha}}$.   
\end{lemma}
For the proof see Appendix \ref{app}.\medskip
\par The $L^1$-continuity of the map $\alpha \mapsto \rho_\alpha$ is sufficient to prove that the entropy map $\alpha \mapsto h(\alpha)$ is also continuous. This is achieved by applying the following lemma (for a proof see for example \cite{AOT}) to Rohlin's formula. 

\begin{lemma} \label{AOT}
Let $\{\rho_n\}$ be a sequence of functions in $L^1(I)$ such that 
\begin{enumerate}
\item $\norm{\rho_n}_\infty \leq K \quad \forall n$, 
\item $\rho_n  \xrightarrow{\;L^1\;} \rho$ for some $\rho \in L^1(I)$ 
\end{enumerate}
Then for any $\psi \in L^1(I)$, 
\begin{equation*}
\int \psi (\rho_n-\rho) \to 0
\end{equation*}
\end{lemma} 

Applying Rohlin's Formula for the entropy, we get for any $\alpha \in [\bar{\alpha}-\varepsilon,\bar{\alpha}+\varepsilon]$
\begin{equation*}
h(\alpha)= \int_{\alpha-1}^{\alpha} \log\frac{1}{(x-\bar{\alpha}+\alpha)^2} \tilde{\rho}_\alpha(x) dx=
2 \int_{\alpha-1}^{\alpha} \abs{\log\abs{x-\bar{\alpha}+\alpha}} \tilde{\rho}_\alpha(x) dx
\end{equation*} 
Consider a sequence $\{\alpha_n\} \to \bar{\alpha}$. Then 
\begin{align*}
&\abs{h(\bar{\alpha})-h(\alpha_n)} \leq 2\int_{\bar{\alpha}-1}^{\bar{\alpha}} \big |\log\abs{x-\bar{\alpha}+\alpha_{n}}\tilde{\rho}_{\alpha_n}(x)-\log\abs{x} \tilde{\rho}_{\bar{\alpha}}(x)\big|dx \leq \\
&\leq 2\left(\int_{\bar{\alpha}-1}^{\bar{\alpha}} \big|\log\abs{x-\bar{\alpha}+\alpha_{n}}\left(\tilde{\rho}_{\alpha_n}(x)-\tilde{\rho}_{\bar{\alpha}}(x)\right)\big|dx +\right.\\
&\left.+\int_{\bar{\alpha}-1}^{\bar{\alpha}} \big|\left(\log\abs{x-\bar{\alpha}+\alpha_{n}}-\log\abs{x}\right) \tilde{\rho}_{\bar{\alpha}}(x)\big|dx\right) 
\end{align*} 
The second integral is bounded by $2 (K_0+1) \int_{\bar{\alpha}-1}^{\bar{\alpha}} \abs{\log\abs{x-\bar{\alpha}+\alpha_{n}}-\log\abs{x}} dx$ and vanishes when $n \to \infty$ because of the continuity of translation in $L^1$. 
If we take $\tilde{\rho}_n=\tilde{\rho}_{\alpha_n}$, $\tilde{\rho}=\tilde{\rho}_{\bar{\alpha}}$, $\psi(x)=\abs{\log\abs{x}}$ in Lemma \ref{AOT}, we find that the first integral also tends to $0$.\qed

\subsection{Behaviour of the density and entropy when $\alpha \to 0$}

In this section we will prove that the entropy has a limit as $\alpha \to 0^+$ and that $\lim_{\alpha \to 0^+} h(\alpha)=0$.\\
The continuity of the entropy on the interval $(0, \sqrt{2}-1]$ followed from the $L^1$-continuity of the densities. The vanishing of the entropy as $\alpha \to 0$ is a consequence of the fact that the densities converge to the Dirac delta at the parabolic fixed point of $T_0$ as $\alpha \to 0$. 
\begin{proposition} \label{delta}
When $\alpha \to 0$, the invariant measures $\tilde{\mu}_{\alpha}$ of the translated maps $A_{\alpha,0}: [-1,0] \to [-1,0]$ converge in the sense of distributions to the Dirac delta in $-1$.
\end{proposition}

From the previous Proposition the vanishing of the entropy follows easily:

\begin{corollary} \label{zero}
Let $h(\alpha)$ be the metric entropy of the map $T_{\alpha}$ with respect to the absolutely continuous invariant probability measure $\mu_{\alpha}$. Then $h(\alpha)\to 0$ as $\alpha \to 0$.
\end{corollary}
\begin{proof}[\textbf{Proof of the Corollary}]
We compute the entropy of the $T_\alpha$ through Rohlin's formula:  
\begin{equation} \label{entropy}
h(\alpha)=2\int_{\alpha-1}^{\alpha} \abs{\log{\abs{x}}}d\mu_{\alpha}
\end{equation}
Observe that $\forall E \subseteq (c_1,0]$, $\mu_\alpha(E)=\frac{1}{C(\alpha)} \nu_{\alpha}(E)\leq \frac{C_0}{C(\alpha)} m(E)$. Therefore if $\rho_\alpha$ is the density of $\mu_{\alpha}$, $\rho_\alpha < \frac{C_0}{C(\alpha)}$ in $(c_1,0]$. 
Given $\varepsilon$, let $c_k$ be such that $\abs{\log\abs{x}} < \varepsilon$ for $x \in [-1,c_k]$, and choose $\alpha$ small such that $\alpha-1 < c_k$, $\mu_{\alpha}([c_k,\alpha])=\tilde{\mu}_{\alpha}([\tilde{c}_k,0])<\varepsilon$ and $\frac{C_0}{C(\alpha)}<\varepsilon$. Then
\begin{multline*}
h(\alpha) \leq \int_{\alpha-1}^{c_k} \abs{\log\abs{x}} d\mu_\alpha+\int_{c_k}^{c_1} \abs{\log{\abs{x}}} d\mu_\alpha+\int_{c_1}^{\alpha} \abs{\log\abs{x}} \rho_{\alpha} dx \leq \\
\leq \abs{\log \abs{c_k}}+ \abs{\log{\frac{1}{3}}} \mu_\alpha([c_k,c_1])+\frac{C_0}{C(\alpha)} \norm{\log{\abs{x}}}_1 \to 0
\end{multline*} 
which concludes the proof. 
\end{proof}

To prove Proposition \ref{delta} we adopt the following strategy: we introduce the jump transformations $G_{\alpha}$ of the maps $T_{\alpha}$ over the cylinder $(2,-)$, whose derivatives are strictly bounded away from $1$ even when $\alpha \to 0$; we can then prove that their densities $\frac{d\nu_{\alpha}}{dx}$ are bounded from above and from below by uniform constants. Using the relation between $\mu_{\alpha}$ and the induced measure $\nu_{\alpha}$, we conclude that for any measurable set $B$ such that $-1 \notin B$, $\tilde{\mu}_{\alpha}(B)=\mu_{\alpha}(B+\alpha) \to 0$ when $\alpha \to 0$. 

\begin{proof}[\textbf{Proof of Proposition \ref{delta}}] 
Given $v_{r+1} \leq \alpha < v_r$ as in equation (\ref{v_r}), and $0 \leq j \leq r$, let
\begin{displaymath}
L_0=I_{\alpha}\setminus (2,-),\quad L_j=[c_{j+1},c_j)=(\underbrace{(2,-),\ldots,(2,-)}_{j})\setminus (\underbrace{(2,-),\ldots,(2,-)}_{j+1})
\end{displaymath}
for $1 \leq j \leq r$. Thus $I_{\alpha}=\bigcup_{0 \leq j \leq r} L_j \pmod{0}$. It is easy to prove by induction that for $r \geq j \geq 1,\; c_j=V_{(2,-)}^{j-1}\left(-\frac{1}{2+\alpha}\right)=-1+\frac{1}{j+\frac{1}{1+\alpha}}$, that is, $-\frac{j}{j+1} < c_j \leq -\frac{j-1}{j}$, while $c_0=\alpha,\,c_{r+1}=\alpha-1$. Let
\begin{displaymath}
G_{\alpha}|_{L_j}=T_{\alpha}^{j+1}|_{L_j}
\end{displaymath}
be the jump transformation associated to the return time $\tau(x)=j+1 \iff x \in L_j$. Observe that $\tau$ is bounded and therefore integrable with respect to $\mu_\alpha$. Then a result of R. Zweimuller (\cite{Zw2}, Theorem 1.1) guarantees that $G_\alpha$ admits an invariant measure $\nu_\alpha \ll \mu_\alpha$ such that for all measurable $E$,
\begin{equation} \label{induced_measure}
\mu_{\alpha} (E) =\frac{1}{C(\alpha)} \left( \sum_{n \geq 0} \nu_{\alpha}\left(\{ \tau >n\} \cap T_{\alpha}^{-n}(E)\right)\right)
\end{equation}  
where $C(\alpha)$ is a suitable normalization constant. Actually from equation (\ref{induced_measure}) it follows that $\nu_{\alpha}(I_{\alpha})=\nu_{\alpha}(\{\tau>0\})\leq C(\alpha) \mu_{\alpha}(I_{\alpha})$ is finite, and so by choosing a suitable $C(\alpha)$ we can take $\nu_{\alpha}(I_{\alpha})=1$. We will prove the following:

\begin{lemma} \label{constants}
\emph{There exists $\tilde{\alpha}>0$ such that for $0 < \alpha<\tilde{\alpha}$, the densities $\psi_{\alpha}$ of $\nu_\alpha$ are bounded from above and from below by constants that do not depend on $\alpha$:} $\exists C_0$ s.t. $C_0^{-1} \leq \psi_{\alpha} \leq C_0$.
\end{lemma} 
\begin{proof}[\textbf{Proof of Lemma \ref{constants}}]
In order to prove that $\psi_{\alpha}$ is bounded from above, we can proceed as in Lemma \ref{Lemma}, and show that $\exists C'$ such that for all $\alpha$, $\forall \varphi \in L^1(I_\alpha)$, $\Var_{I_{\alpha}}P_{G_\alpha}^n \varphi < C'$. Since the outline of the proof is very similar to that of Lemma \ref{Lemma}, we will only list the passages where the estimates are different, and emphasize how in this case all the constants can be chosen uniform in $\alpha$.\\ 
The cylinders of rank $1$ for $G_\alpha$ are of the form
\begin{displaymath}
I_{j}^{k,\varepsilon}=(j,k,\varepsilon)\doteqdot(\underbrace{(2,-),\ldots,(2,-)}_{j},(k,\varepsilon)),\quad 0 \leq j \leq r,
\end{displaymath} 
so they are also cylinders for $T_\alpha$, although of different rank. On $I_{j}^{k,\varepsilon},\; j \geq 1$ we have 
\begin{align*}
&\abs{\frac{1}{G_\alpha'(x)}}=(T_{\alpha}^j(x) \cdots T_\alpha(x) x)^2 \leq \lambda_{j}^{k}= \frac{4}{(j+2)^2 (k-1)^2}\leq \frac{1}{(k-1)^2} \leq \frac{1}{4},\\
& \abs{\frac{1}{G_\alpha'(x)}} \geq \frac{1}{9 j^2 (k+1)^2},
\end{align*}
while on $I_{0}^{k,\varepsilon}$, $\abs{\frac{1}{G_{\alpha}'(x)}}=x^2 <\frac{1}{(k-1)^2} \leq \frac{1}{4}$, and so $\lambda=\sup \abs{\frac{1}{G_\alpha'}} <\frac{1}{4}$ for all $\alpha$.\\
Letting $\mathcal{Q}=\bigcup_{j=0}^{r}\{I_{j}^{k,\varepsilon}\}$, and $\lambda_{\eta}^{(n)}=\sup_{I_{\eta}^{(n)}} \abs{\frac{1}{(G_{\alpha}^n)'}}$, we can obtain the analogue of equation (\ref{var_P}) for the maps $G_{\alpha}$:
\begin{equation*} 
\Var \left(P_{G_{\alpha}}^n \varphi \right)(x) \leq \lambda^n \Var_{I_\alpha} \varphi + \sum_{I_{\eta}^{(n)} \in \mathcal{Q}^{(n)}} (n +2) \lambda_{\eta}^{(n)} \sup_{I_{\eta}^{(n)}}\abs{\varphi},  
\end{equation*} 
and similarly to (\ref{k}), we can choose $h$ such that $\sum_{i \geq h} \frac{1}{i^2} \leq \frac{\lambda^n}{2^{4n-2}}$,
and the set of intervals with bounded itineraries
\begin{multline*} 
G(n)=\{ I_{\eta}^{(n)}=((j_0,k_0,\varepsilon_0),\ldots,(j_{n-1},k_{n-1},\varepsilon_{n-1})) \in \mathcal{Q}^n \;| \\ \max(j_0,\ldots,j_{n-1}) \leq h,\;\max(k_0,\ldots,k_{n-1}) \leq h \}
\end{multline*}
Again we can define a function $\sigma: G(n) \to G(n)$ that maps every cylinder $I_{\eta}^{(n)}=((j_0,k_0,\varepsilon_0),\ldots,(j_{n-1},k_{n-1},\varepsilon_{n-1}))$ to $I_{\xi}^{(n)}=((j_{0}',k_{0}',\varepsilon_{0}'),\ldots,(j_{n-1}',k_{n-1}',\varepsilon_{n-1}'))$ as follows:
\begin{enumerate}
\renewcommand{\theenumi}{\alph{enumi}}
\item If $(j_i,k_i,\varepsilon_i)=(j,j_{\min},+),\; j<r-1$ for some $i$, then $(j_{i}',k_{i}',\varepsilon_{i}')=(j,j_{\min}+1,+)$;
\item If for some $i$, $(j_i,k_i,\varepsilon_i)=(r,k,\varepsilon)$ with $(k,\varepsilon) \neq (j_{\min},+),(j_{\min}+1,+)$, then $(j_{i}',k_{i}',\varepsilon_{i}')=(r-1,k,\varepsilon)$;
\item If $(j_i,k_i,\varepsilon_i) \in \{(r,j_{\min},+),(r,j_{\min}+1,+),(r-1,j_{\min},+)\}$, then\\ $(j_{i}',k_{i}',\varepsilon_{i}')=(r-1,j_{\min}+1,+)$;
\item Otherwise, $(j_{i}',k_{i}',\varepsilon_{i}')=(j_i,k_i,\varepsilon_i)$.
\end{enumerate}
With this definition, the cylinders in $\sigma(G(n))$ are all full, because as we have seen in equation (\ref{surj}), for $0 \leq i \leq r-1$,
\begin{equation*}
T_{\alpha}^i(\underbrace{(2,-),\ldots,(2,-)}_{i})\supseteq \left[-\frac{1}{2+\alpha},\alpha\right)\supseteq \bigcup\limits_{(k,\varepsilon), k \geq 3} I_{k}^{\varepsilon} 
\end{equation*}
Then $\forall I_{\xi}^{(n)} \in \sigma(G(n))$, 
\begin{equation*}
1 \leq m(I_{\xi}^{(n)}) \sup_{\sigma(G(n))} \abs{(G_{\alpha}^n)'}= m(I_{\xi}^{(n)}) (9h^4)^n \Rightarrow m(I_{\xi}^{(n)}) \geq \frac{1}{\delta_n}=\frac{1}{(9h^4)^n},
\end{equation*}
which doesn't depend on $\alpha$. Again we need to estimate the supremum over $I_{\xi}^{(n)} \in \sigma(G(n))$ of the sums $\sum\limits_{\sigma(I_{\eta}^{(n)})=I_{\xi}^{(n)}}\lambda_{\eta}^{(n)}$, each of which
is the product of $n$ terms $\lambda'_i$, that correspond to one of the cases a),\,b),\,c)\,d) we listed previously: 
\begin{itemize}
\item In the case a), $\lambda'_i=\lambda_{j}^{j_{\min}}+\lambda_{j}^{j_{\min}+1} \leq \frac{1}{(j_{\min}-1)^2}+\frac{1}{j_{\min}^2}\leq \frac{1}{2}$ (observe that for $\alpha<\sqrt{2}-1$, $j_{\min} \geq 3$). 
\item In the case b), $\lambda'_i=\lambda_{r}^{k}+\lambda_{r-1}^{k} \leq \frac{4}{(k-1)^2}\left(\frac{1}{(r+1)^2}+\frac{1}{(r+2)^2}\right) \leq \frac{1}{2}$ when $\alpha<v_2$; 
\item In the case c), $\lambda'_i=\lambda_{r}^{j_{\min}}+\lambda_{r-1}^{j_{\min}}+\lambda_{r}^{j_{\min}+1}+\lambda_{r-1}^{j_{\min}+1}\\
<4\left(\frac{1}{(j_{\min}-1)^2}+\frac{1}{j_{\min}^2}\right)\left(\frac{1}{(r+1)^2}+\frac{1}{(r+2)^2}\right)<\frac{1}{2}$ for $\alpha<v_2$;
\item In the case d), $\lambda'_i< \lambda=\frac{1}{4}$. 
\end{itemize}
Then $\lambda'\leq \tilde{\lambda}=\frac{1}{2}$, and as in equation (\ref{tilde}), we find for $\alpha< v_2$,
\begin{equation*}
\sum_{I_{\eta}^{(n)} \in G(n)} \left( \lambda_{\eta}^{(n)} \sup_{I_{\eta}^{(n)}}\abs{\varphi} \right) \leq \lambda' \sum_{I_{\xi}^{(n)} \in \sigma(G(n))} \sup_{\overline{I}_{\xi}^{(n)}} \abs{\varphi} 
\leq \tilde{\lambda}^n \Var \varphi + \frac{\tilde{\lambda}^n}{\delta_n} \norm{\varphi}_1  
\end{equation*} 
For the sum over intervals with unbounded itineraries we proceed in a similar way to (\ref{9}):
\begin{multline*}
\sum_{I_{\eta}^{(n)} \notin G(n)} \lambda_{\eta}^{(n)} \leq \sum_{i\geq h} \sum_{l=0}^{n-1}\left( \sum_{\substack{j_{l}(\eta)=i+1=\\ \max (j_0(\eta),\ldots,j_{n-1}(\eta))}} \lambda_{j_0(\eta)}^{k_0(\eta)} \cdots \lambda_{j_{n-1}(\eta)}^{k_{n-1}(\eta)}+ \right.\\  +\left. \sum_{\substack{k_{l}(\eta)=i=\\\max (k_0(\eta),\ldots,k_{n-1}(\eta))}} \lambda_{j_0(\eta)}^{k_0(\eta)} \cdots \lambda_{j_{n-1}(\eta)}^{k_{n-1}(\eta)}  \right) \leq \sum_{i\geq h} \sum_{l=0}^{n-1} \frac{4}{i^2} \left(2\sum_{I_{j}^{k,\varepsilon} \in \mathcal{Q}} \lambda_{\eta}^{(n)}\right)^{n-1}
\end{multline*}
(This expression is redundant, but sufficient for our purpose.)
Observe that 
\begin{align*}
&\sum_{I_{j}^{k,\varepsilon} \in \mathcal{Q}} \lambda_{\eta}^{(n)} \leq 2 \sum_{j=0}^\infty \sum_{k=3}^{\infty} \frac{4}{(j+2)^2 (k-1)^2} \leq 8\left(\sum_2^\infty \frac{1}{k^2}\right)^2 \leq 8, \text{ and so }\\
&\sum_{I_{\eta}^{(n)} \notin G(n)} \lambda_{\eta}^{(n)} \leq \sum_{i\geq h} \frac{4}{i^2} n 2^{4n-4} \leq n \lambda^n
\end{align*}
Then we can prove relation (\ref{1a}) and complete our argument exactly like in Lemma 3. Notice that all the constants involved are uniform in $\alpha$.\bigskip

\par To prove that the densities $\psi_\alpha$ of $\nu_{\alpha}$ are uniformly bounded from below, we use a bounded distortion argument. We follow the same outline as in Appendix \ref{bounded_distortion}, but with the advantage that in this case the derivatives are uniformly bounded from above.\\ 
Since $T_\alpha$ satisfies \emph{Adler's condition} $\abs{\frac{T_{\alpha}''}{(T_{\alpha}')^2}}<K$ (here $K=2$), then there exists $K'$ independent of $n$ and of $\alpha$ such that $\forall n>0,\, \abs{\frac{(T_{\alpha}^n)''}{((T_{\alpha}^n)')^2}}<K'$ (see \cite{Zw1}, Lemma 10). Then $\forall x,y$ belonging to the same cylinder $I_{j}^{k,\varepsilon}$ of rank $1$ of $G_\alpha$,
\begin{multline*}
\abs{\frac{G_{\alpha}'(x)}{G_{\alpha}'(y)}-1}= \abs{G_{\alpha}''(\xi)} \abs{\frac{x-y}{G_\alpha'(y)}}= \abs{\frac{G_{\alpha}''(\xi)(G_\alpha(x)-G_\alpha(y))}{G_\alpha'(y) G_\alpha'(\eta)}} \leq \\
\leq 36^2 \abs{\frac{G_\alpha''(\xi)}{(G_\alpha'(\xi))^2}} \abs{G_\alpha(x)-G_{\alpha}(y)} \leq K'' \abs{G_\alpha(x)-G_{\alpha}(y)}, \quad \text{ and} 
\end{multline*}
\begin{multline*} 
\log \abs{\frac{(G_{\alpha}^n)'(y)}{(G_{\alpha}^n)'(x)}} \leq \sum_{i=0}^{n-1} \abs{\frac{G_{\alpha}'(G_{\alpha}^i(y))}{G_{\alpha}'(G_{\alpha}^i(x))}-1}  
\leq K'' \sum_{i=0}^{n-1} \abs{G_{\alpha}^{i+1}(y)-G_{\alpha}^{i+1}(x)} \leq \\
 \leq K'' \sum_{i=1}^n \left(\frac{1}{4}\right)^{n-i} \abs{G_{\alpha}^{n}(y)-G_{\alpha}^{n}(x)} \leq 
K'' \sum_{i=0}^\infty \left(\frac{1}{4}\right)^i \;\Rightarrow\; \abs{\frac{(G_{\alpha}^n)'(y)}{(G_{\alpha}^n)'(x)}} \leq C_1
\end{multline*}
where we remark that $C_1$ does not depend on $\alpha$. Letting $W_{\eta} : G_{\alpha}^n(I_{\eta}^{(n)}) \to I_{\eta}^{(n)}$ be the local inverses of $G_{\alpha}^n$, for every full cylinder $I_{\eta}^{(n)} \in \mathcal{P}^{(n)}$ and for every measurable set $B$,
\begin{multline} \label{S_n} 
\frac{m(B)}{m(I_{\alpha})}=\frac{\int\limits_{W_\eta(B)}\abs{(G_{\alpha}^n)'(y)}dy}{\int\limits_{I_{\eta}^{(n)}}\abs{(G_{\alpha}^n)'(x)}dx} \leq \frac{m(W_\eta(B)) \sup\limits_{y \in I_{\eta}^{(n)}} \abs{(G_{\alpha}^n)'(y)}}{m(I_{\eta}^{(n)}) \inf\limits_{x \in I_{\eta}^{(n)}} \abs{(G_{\alpha}^n)'(x)}} \leq C_1 \frac{m(W_\eta(B))}{m(I_{\eta}^{(n)})} \\
\Rightarrow m(W_{\eta}(B)) \geq m(B) \frac{m(I_{\eta}^{(n)})}{C_1} 
\end{multline}
Finally, we can show that the measure of the union $S_n$ of all full cylinders of rank $n$ is strictly greater than $0$. In fact we have the following characterization:
$I_{\eta}^{(n)} \in \mathcal{Q}^{(n)}$ is not full $\Rightarrow$ it has an initial segment of the orbit (with respect to $G_{\alpha}$) of one of the endpoints $\alpha$ and $\alpha-1$ as its final segment. That is, if $\alpha=(a_1,a_2,a_3,\ldots)$ and $\alpha-1=(b_1,b_2,b_3,\ldots)$, then there exists $1 \leq k \leq n$ such that $I_{\eta}^{(n)}=(\omega_1,\ldots,\omega_{n-k},a_1,\ldots,a_k)$ or $I_{\eta}^{(n)}=(\omega_1,\ldots,\omega_{n-k},b_1,\ldots,b_k)$. To prove this, observe that if $I_{\eta}^{(n)}$ doesn't contain any initial segment of $(a_1,a_2,a_3,\ldots)$ or $(b_1,b_2,b_3,\ldots)$, it is clearly full, and if every such segment $(a_1,\ldots,a_k)$ or $(b_1,\ldots,b_k)$ is followed by $\omega_{k+1} \neq a_{k+1}$ or $b_{k+1}$ respectively, then it is either full or empty because $G_{\alpha}^n$ is monotone on each cylinder. Then 
\begin{align*} 
& \nu_{\alpha}(S_n^\mathcal{C})\leq \nu_{\alpha}\left(\bigcup_{k=1}^{n} G_{\alpha}^{-(n-k)} (a_1,\ldots,a_k)\right)+\nu_{\alpha}\left(\bigcup_{k=1}^{n} G_{\alpha}^{-(n-k)} (b_1,\ldots,b_k)\right) \leq \\
& \leq \sum_{k=1}^n \left(\nu_{\alpha}(a_1,\dots,a_k)+\nu_{\alpha}(b_1,\ldots,b_k)\right)
\end{align*} 
since $\nu_{\alpha}$ is $G_{\alpha}$-invariant. We have already shown that $\nu_\alpha$ is bounded from below, and so $\nu_{\alpha}(a_1,\ldots,a_k)+\nu_{\alpha}(b_1,\ldots,b_k) < C' (m(a_1,\ldots,a_k)+m(b_1,\ldots,b_k))$.\\
In order to prove that $\nu_\alpha(S_n^C)<1$, we take advantage of the fact that the cylinders containing the endpoints become arbitrarily small when $\alpha$ approaches $0$. Recall that $(a_1)=(j_{\min})=\left(\left[\frac{1}{\alpha}+1-\alpha\right]\right)$, and consequently $\inf_{(a_1)} \abs{G_{\alpha}'(x)} \geq (j_{\min}-1)^2$, and since $\inf_{I_{\alpha}} \abs{G_{\alpha}'(x)} \geq 4$, from Lagrange's theorem we get 
\begin{equation*}
m(a_1,\ldots,a_k) \leq \frac{1}{4^{k-1} (j_{\min}-1)^2}
\end{equation*}  
Since $j_{\min} \to \infty$ as $\alpha \to 0$, we can choose $\tilde{\alpha}$ such that $\forall \alpha<\tilde{\alpha},\; m(a_1,\ldots,a_k) \leq \frac{1}{4^{k} C'}$. Similarly, recall from Remark \ref{nu} that for $v_{r+1} \leq \alpha < v_r$, where $v_{r}=\frac{-1+\sqrt{1+4/r}}{2}$, we have $(b_1)=(\underbrace{(2,-),\ldots,(2,-)}_{r},(k,\varepsilon))$ for some $k \geq 3$, and recalling that $T_{\alpha}^i(\alpha-1)=\frac{(i+1)\alpha-1}{1-i\alpha}$, we find
\begin{equation*}
\inf_{(b_1)} \abs{G_{\alpha}'(x)}=\inf_{(b_1)} \prod_{i=0}^r\frac{1}{(T_{\alpha}^i(x))^2} \geq \frac{1}{(k-1)^2}\prod_{i=0}^{r-1} \frac{(1-i\alpha)^2}{(1-(i+1)\alpha)^2} \geq \frac{4}{(1-r\alpha)^2} 
\end{equation*} 
But $\alpha \geq v_{r+1} \Rightarrow \frac{1}{r+1}<\alpha^2+\alpha \Rightarrow 1-r\alpha< \frac{\alpha(2+\alpha)}{1+\alpha} < 3 \alpha$, and so by taking $\tilde{\alpha}$ small enough we can ensure that $\forall \alpha<\tilde{\alpha},\;\inf_{(b_1)} \abs{G_{\alpha}'(x)} \geq 4C'$ and consequently $m(b_1,\ldots,b_k) \leq \frac{1}{4^{k} C'} \;\forall k$.\\   
Then for $\alpha$ small enough, $m(S_n^\mathcal{C})\leq \frac{2}{C'} \sum_{k=1}^n \frac{1}{4^k} \leq \frac{2}{3 C'} \Rightarrow \nu_{\alpha}(S_n) \leq \frac{2}{3} \Rightarrow \nu_{\alpha}(S_n) \geq \frac{1}{3} \Rightarrow m(S_n) > \frac{\nu_{\alpha}(S_n)}{C'} \geq \frac{1}{3 C'}$. Taking the sum over all full cylinders $I_{\eta}^{(n)}$ in (\ref{S_n}), we find that for all measurable $B \subseteq I_{\alpha}$,
\begin{equation*}  
m(G_{\alpha}^{-n}(B)) \geq m(G_{\alpha}^{-n}(B) \cap S_n) \geq \frac{m(B) m(S_n)}{C_1} \geq \frac{m(B)}{3 C_1 C'}
\end{equation*}
Now recall that the density of $\nu_{\alpha}$ is equal almost everywhere to the limit of the Cesaro sums $\lim_{n \to \infty} \frac{1}{n} \sum_{i=0}^{n-1} P_{G_{\alpha}}^{\;i} 1$, and so for $\alpha < \tilde{\alpha}$,
\begin{equation*}
\nu_\alpha(B)=\lim_{n \to \infty} \frac{1}{n} \sum_{i=0}^{n-1} \int\limits_B P_{G_{\alpha}}^{\;i} 1 dx= 
\lim_{n \to \infty} \frac{1}{n} \sum_{i=0}^{n-1} m(G_{\alpha}^{-n}(B)) 
\end{equation*}
and consequently we have $\nu_\alpha(B) \geq  \frac{m(B)}{3 C_1 C'}$ $\forall B$.
\end{proof}

We can finally conclude the proof of Proposition \ref{delta}. The following properties hold:
\begin{itemize}
\item \emph{$C(\alpha) \to \infty$ when $\alpha \to 0$.} In fact when $\alpha$ is small, $C_0^{-1} \leq \frac{d\nu_\alpha}{dm} \leq C_0$ for some $C_0$, and
\begin{multline} \label{Viana}
1= \sum_{k=0}^r \mu_{\alpha} (L_k)= \sum_{k=0}^r \frac{1}{C(\alpha)} \sum_{n\geq k} \nu_\alpha(L_{n})=\frac{1}{C(\alpha)} \sum_{k=0}^r \nu_\alpha ([\alpha-1,c_k]) \geq \\
\geq \frac{1}{C_0 C(\alpha)} \sum_{k=0}^r m([\alpha-1,c_k]) \geq \frac{1}{C_0 C(\alpha)}\left ( \sum_{k=0}^r \frac{1}{k+1} -(r+1)\alpha \right) \geq \\
\geq \frac{1}{C_0 C(\alpha)}\left(\log\left(\frac{1}{\alpha^2+\alpha}\right)-1\right)
\end{multline}
since $r \leq \frac{1}{\alpha^2 +\alpha}\leq r+1$. Therefore the normalization constant $C(\alpha) \geq \frac{1}{C_0} \left( \log \left(\frac{1}{\alpha^2+\alpha}\right) -1\right) \to \infty$ when $\alpha \to 0$.  	
\item Finally, \emph{$\forall L_k,\; k \geq 0$ finite, $\mu_{\alpha}(L_k)\to 0$ when $\alpha \to 0$.}
In fact we have
\begin{equation*}
\mu_{\alpha}(L_k) =\sum_{j \geq k} \frac{\nu_{\alpha}(L_j)}{C(\alpha)} \leq \frac{C_0}{C(\alpha)} \sum_{j\geq k} m(L_j) \leq \frac{C_0}{C(\alpha)} \abs{1-c_k}\leq \frac{C_0}{C(\alpha)} \frac{1}{k} \to 0 
\footnote{The reader might be wondering whether the estimates of the densities of $T_\alpha$ in Paragraph \ref{UBV} and the continuity of the entropy might be derived directly from Lemma \ref{constants}. This would follow from equation (\ref{Viana}) if we had a suitable lower bound for $\frac{d\nu_{\alpha}}{dm}$ when $\alpha$ varies; but we haven't been able to provide such a bound except for small $\alpha$. As for the continuity of the entropy, the fact that $h(T_\alpha)$ and $h(G_\alpha)$ are related by the Generalized Abramov Formula \cite{Zw2} suggests that proving the continuity of $h(G_\alpha)$ might be a valid alternative approach; however, taking expansivity into account, we believe that the estimates necessary to prove $L^1$-continuity of the invariant densities of $G_\alpha$ (as in Appendix \ref{app}) would be far more taxing than for $T_\alpha$.} 
\end{equation*}
\end{itemize}
Consider now the translated versions $A_{\alpha,0}$ of the $T_\alpha$ with respect to $\bar{\alpha}=0$, and let $\tilde{c}_j=c_j-\alpha$ be the translated versions of the $c_j$ (we omit the dependence on $\alpha$ for simplicity). Then we have $\tilde{\mu}_\alpha((\tilde{c}_k,0])\to 0$ for all finite $k$. Let $f \in C^{\infty}([-1,0])$ be a test function: we want to show that $\forall \varepsilon>0$, $\exists \alpha'$ such that $\forall \alpha \leq \alpha',\; \abs{\int_{-1}^{0} f(x) d\tilde{\mu}_\alpha-f(-1)}<\varepsilon$. Since $f$ is uniformly continuous, $\exists \delta$ such that $\forall \abs{x-1}<\delta$, $\abs{f(x)-f(-1)}<\varepsilon$. Choose $k$ so that $\tilde{c}_k < -1+\delta$. Then for all $\alpha$ such that $\tilde{\mu}_{\alpha}((\tilde{c}_k,0])< \varepsilon$,  
\begin{align*}
&\abs{\int_{-1}^{0} (f(x)-f(-1)) d\tilde{\mu}_\alpha} \leq \int_{-1}^{\tilde{c}_k} \abs{f(x)-f(-1)} d\tilde{\mu}_\alpha + \int_{\tilde{c}_k}^{0} \abs{f(x)}d\tilde{\mu}_\alpha+\\
&+\int_{\tilde{c}_k}^{0} \abs{f(-1)} d\tilde{\mu}_\alpha\leq \varepsilon+\varepsilon (\norm{f}_\infty+\abs{f(-1)}) \qedhere
\end{align*}  
\end{proof}

\begin{figure}[p] 
\begin{center}
\includegraphics[width=\textwidth, height=0.6\textwidth, angle=0]
{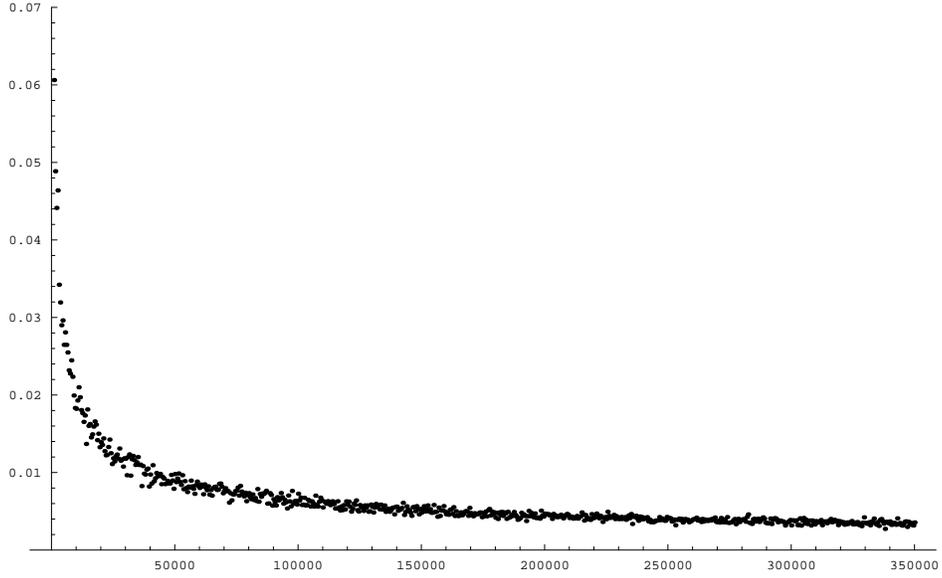}
\caption{The dependence on $n$ of the standard deviation of the normally
distributed $h (\frac{1}{2} , n, x_k)$ where $n$ ranges from $500$ to $350000$ and
$N=100$.}
\label{precision}
\end{center}
\end{figure}

\begin{figure}[p] 
\begin{center}
\includegraphics[width=\textwidth, height=0.6\textwidth, angle=0]
{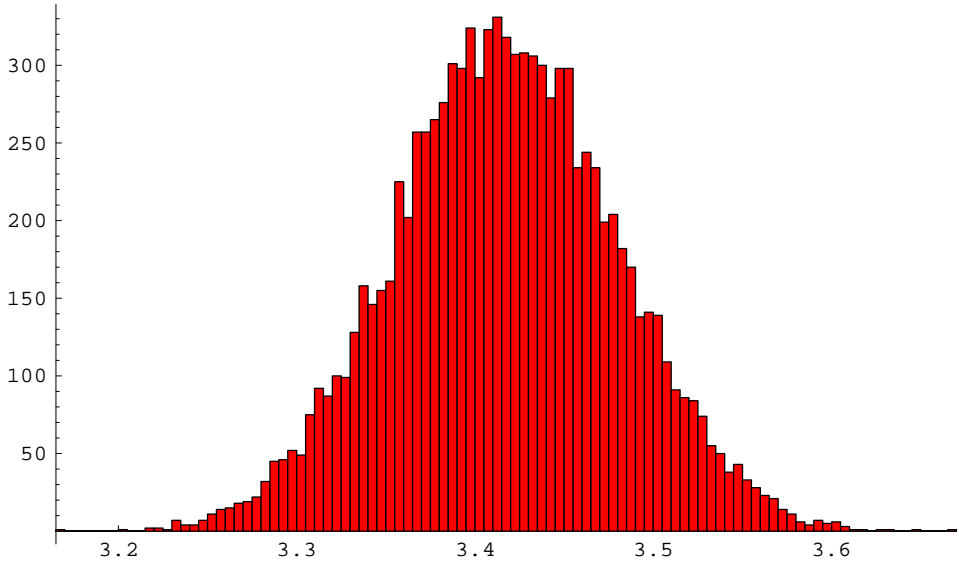}
\caption{The
distribution of $h(\frac{1}{2} , 1000, x_k)$ for $10000$ random initial conditions.
The average $h (\frac{1}{2}, n, N)=3.41711$ must be compared to the exact value
$h(T_{\frac{1}{2}})=\frac{\pi^2}{6\log G}=3.418315971\ldots$}
\label{dependence}
\end{center}
\end{figure}

\begin{figure}[p]
\begin{center}
\includegraphics[width=\textwidth, height=0.6\textwidth, angle=0]
{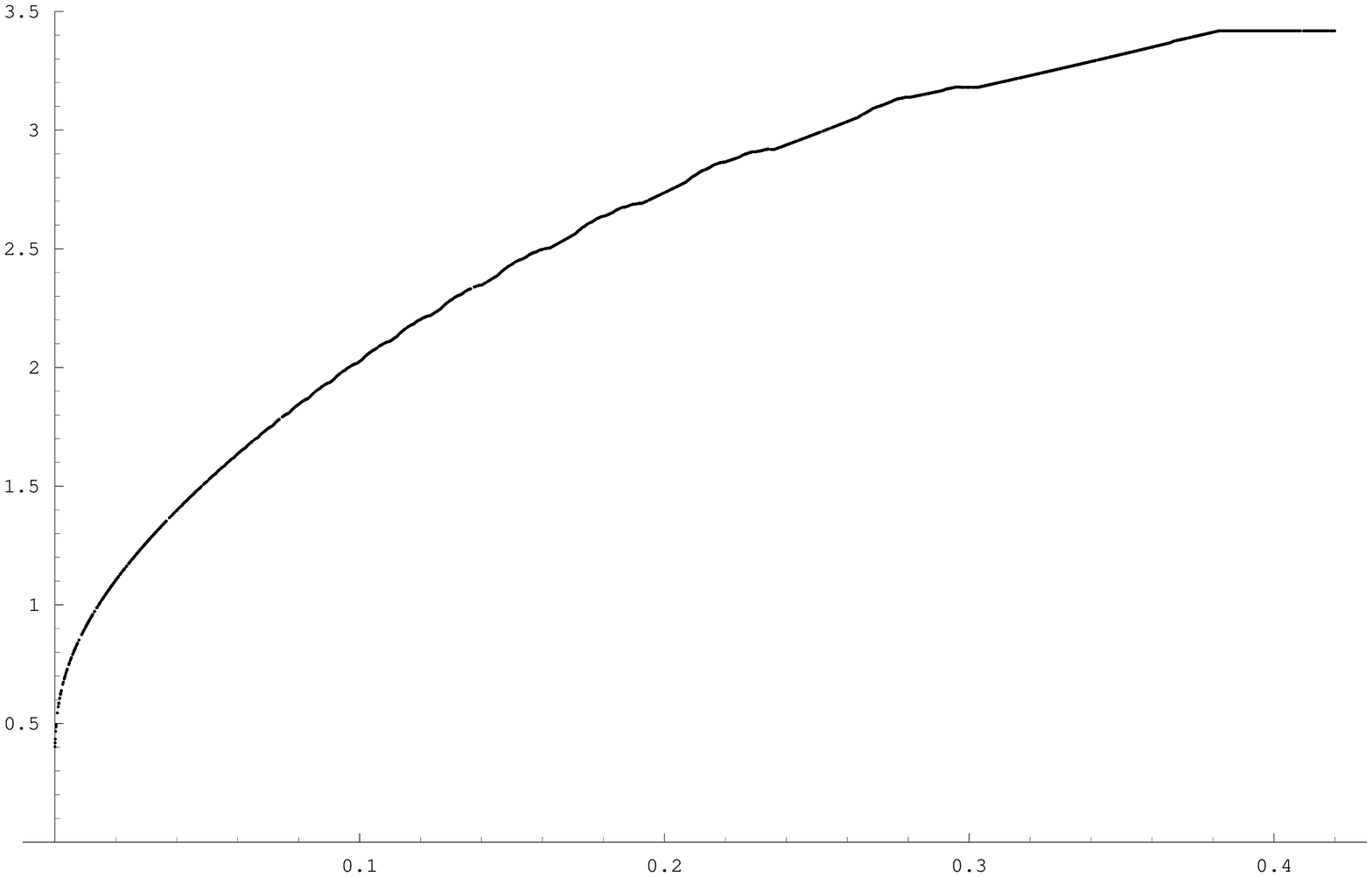}
\caption{The entropy of the map $T_{\alpha}$ at $4080$ uniformly distributed values of $\alpha$ from $0$ to $0.42$. The estimated error is less than $2\cdot 10^{-4}$.}
\label{figure3}
\end{center}
\end{figure}

\begin{figure}[p] 
\begin{center}
\includegraphics[width=\textwidth, height=0.6\textwidth, angle=0]{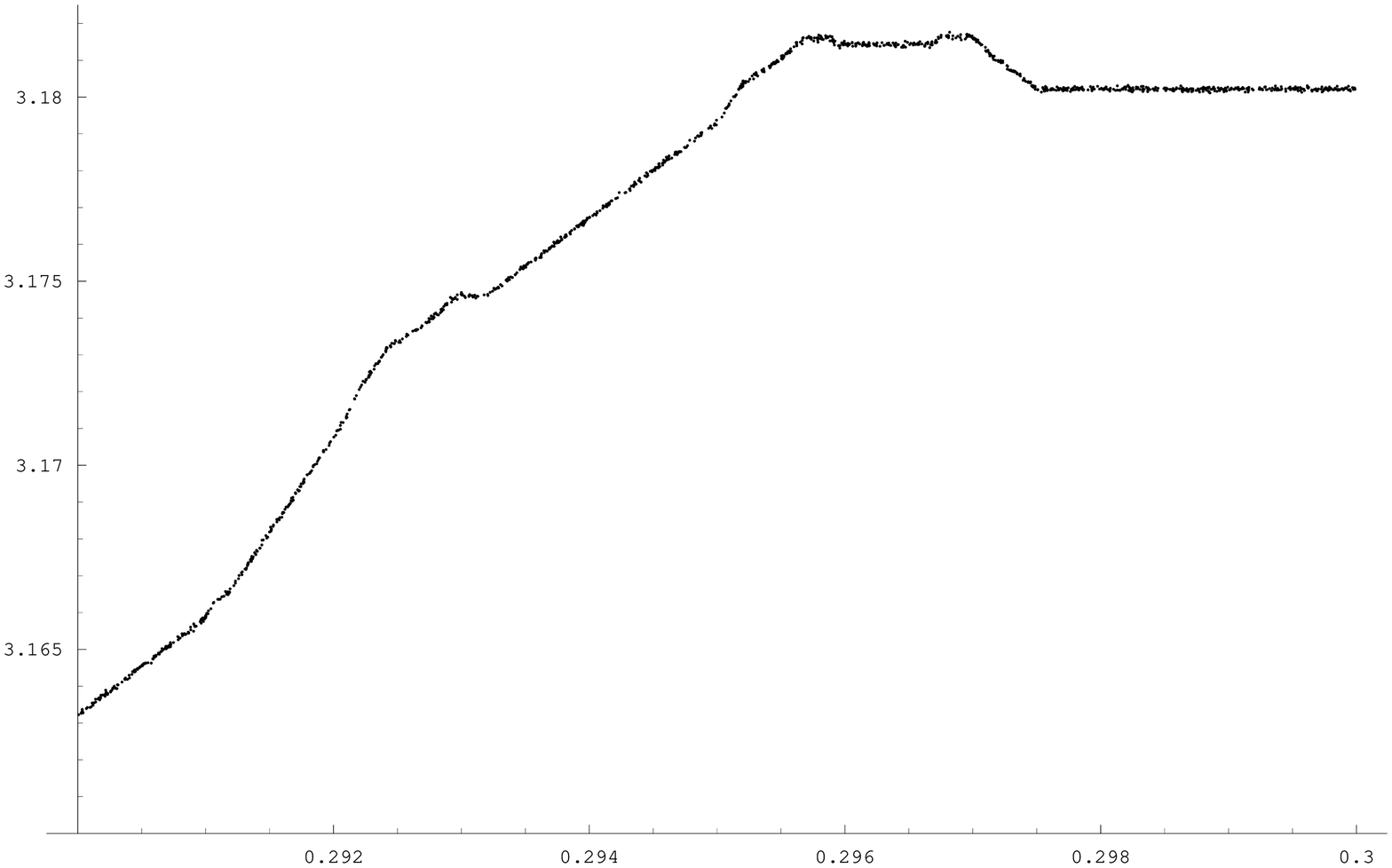}
\caption{The entropy of the map $T_{\alpha}$ at $1314$ uniformly distributed values of $\alpha$ from $0.29$ to $0.30$. The estimated error is less than $1\cdot 10^{-4}$.}
\label{figure4}
\end{center}
\end{figure}

\begin{figure}[p] 
\begin{center}
\includegraphics[width=\textwidth, height=0.6\textwidth, angle=0]{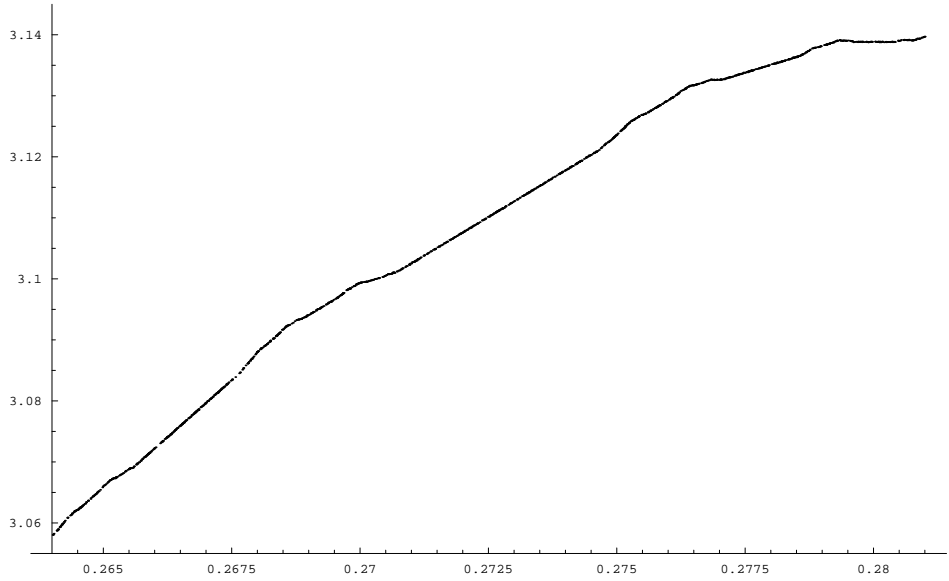}
\caption{The entropy of the map $T_{\alpha}$ at $1600$ uniformly distributed values of $\alpha$ from $0.265$ to $0.281$. The estimated error is less than $1.5\cdot 10^{-4}$.}  
\label{figure5}
\end{center}
\end{figure}

\begin{figure}[p] 
\begin{center}
\includegraphics[width=\textwidth, height=0.6\textwidth, angle=0]{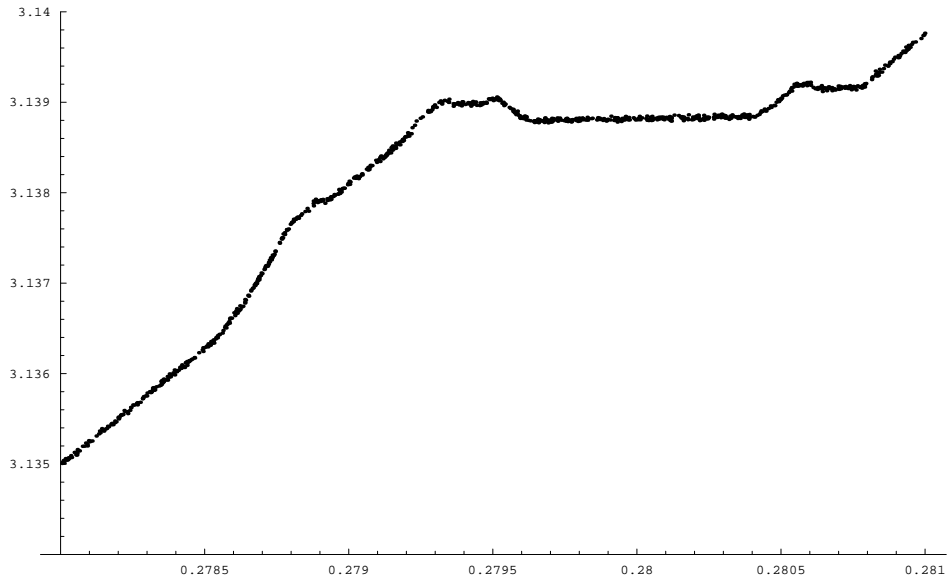}
\caption{The entropy of the map $T_\alpha$ at $989$ 
uniformly 
distributed values of $\alpha$ from $0.278$ to $0.281$. 
The estimated error is less than $4\cdot 10^{-5}$.}
\label{figure6}
\end{center}
\end{figure}

\begin{figure}[htb] 
\begin{center}
\includegraphics[width=\textwidth, height=0.6\textwidth, angle=0]{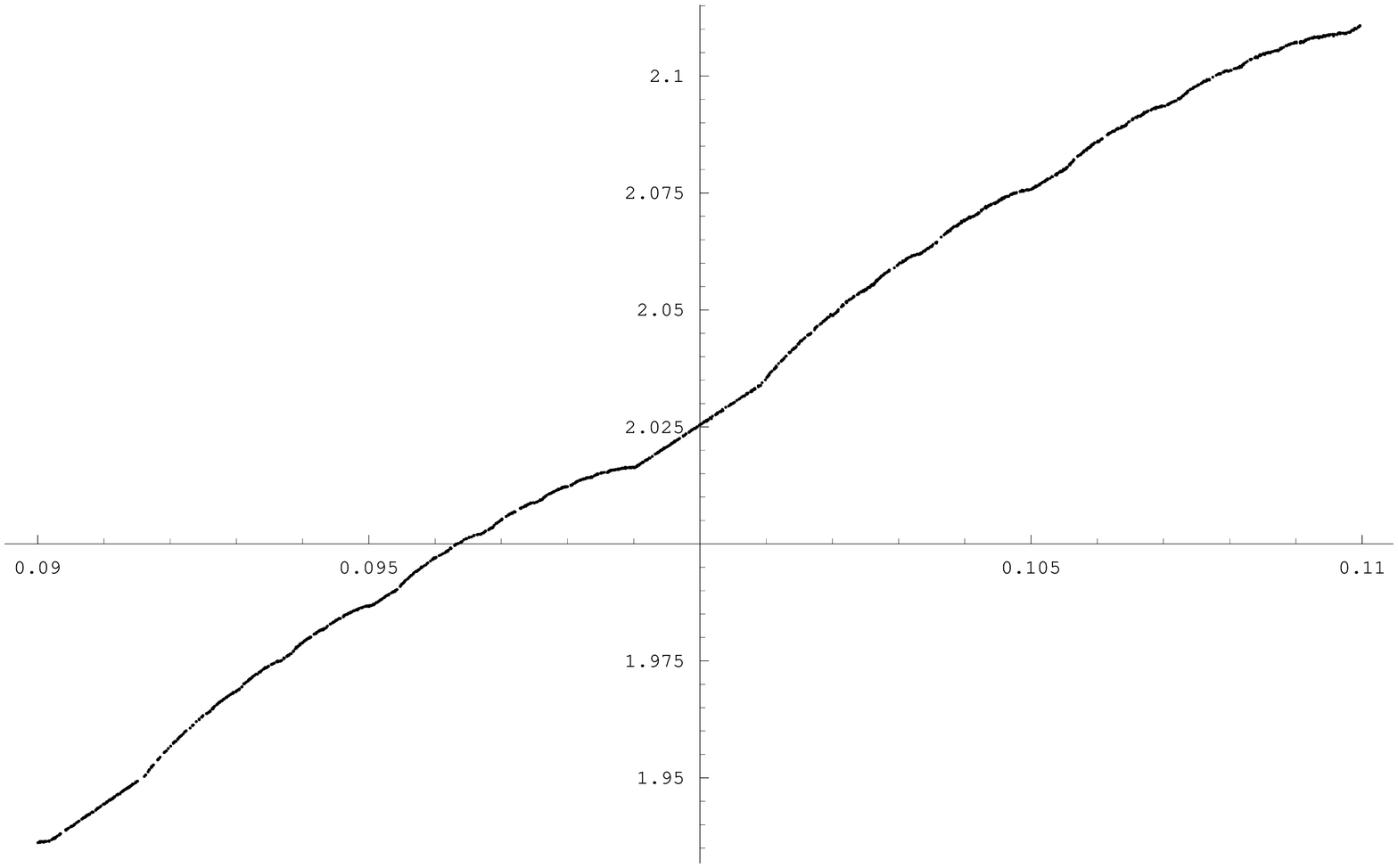}
\caption{The entropy of the map $T_{\alpha}$ at $1799$ uniformly distributed values of $\alpha$ from $0.09$ to $0.11$. The estimated error is less than $2.5\cdot 10^{-4}$.}
\label{figure7}
\end{center}
\end{figure}

\section{Numerical results} \label{numerical_simulations}

In this  section we collect our numerical results on the entropy of Japanese
continued fractions. We already know that the function $\alpha \to h(T_\alpha)$
is con\-ti\-nuous in $(0,1]$ and that in the case $\alpha \geq \sqrt{2}-1$ the entropy
has been computed exactly by Nakada \cite{Nak} and Marmi, Moussa, Cassa
\cite{MCM}. For values of $\alpha$ in the interval $(0,\sqrt{2}-1 ]$ we
have numerically computed the entropy of the maps applying Birkhoff's
ergodic theorem and replacing the integral $h(T_\alpha)= -2 \int_{\alpha
-1}^{\alpha} \log \abs{x} \rho_\alpha(x) dx$ in Rohlin's formula with the Birkhoff
averages
\begin{equation*}
h (\alpha , n, x) = -\frac{2}{n} \sum_{j=0}^{n-1}  \log \abs{T_\alpha^j (x)}
\end{equation*}
which converge to $h(T_\alpha)$ for almost all choices of $x\in (\alpha
-1,\alpha)$. In order to get rid of the dependence on  the choice of an initial
condition we have computed $h (\alpha , n, x_k)$ for a large number $N$ of
uniformly distributed values of  $x_k\in (\alpha -1,\alpha)$, $k=1,\ldots N$,
and we have taken the average on all the results:
\begin{equation*}
h (\alpha , n, N) = \frac{1}{N} \sum_{k=1}^{N}h (\alpha , n, x_k)\; .
\end{equation*}
Unsurprisingly, it turns out that the values $h (\alpha , n, x_k)$ are normally distributed
around their average $h (\alpha , n, N) $ (see Figure \ref{precision}). We have also computed the standard deviations of the normal distributions for values of $n$ from $500$ to $350000$
(see Figure \ref{dependence}): a least squares fit suggests that they decay as $1/\sqrt{n}$
(we refer to A. Broise \cite{Br} for a general treatment of Central Limit Theorems that may apply also to our maps).\\
In Figure \ref{figure3} we see a graph of $h (\alpha , 10^4, N)$ at $4080$ uniformly
distributed random  values of $\alpha$ in the interval $(0,\sqrt{2}-1)$: the
values of $N$ range from $10^5$ to $4\cdot 10^5$ increasing as $\alpha$
decreases so as to keep the standard deviation approximately constant. The estimated error for the entropy is less than $2 \cdot 10^{-4}$.\\ 
As $\alpha \rightarrow 0$ the entropy decreases (although non monotonically, see
below) and the graph exhibits a quite rich self--similar structure that we have
just started to investigate: for example the entropy seems to be independent of $\alpha$ as $\alpha$ varies in the intervals whose endpoints have Gauss continued fraction expansions of the form $[0,\overline{n,n-1,1}]$ and $[0,\overline{n}]$ respectively\footnote{We recall that 
\begin{displaymath}
[0;n_1,n_2,\ldots]=\cfrac{1}{n_1+\cfrac{1}{n_2+\cdots}}
\end{displaymath}
and that the overline indicates the periodic part of the expansion.
}, and to depend linearly on $\alpha$ in the intervals $([0,\overline{n}],[0,\overline{n-1,1}])$. Compare with Figure \ref{figure5}, where $h (\alpha , 10^4,
4\cdot 10^5)$ is computed at $1600$ values of $\alpha\in (0.264, 0.281)$ and
with Figure \ref{figure7} where $h (\alpha , 10^4, 2\cdot 10^5)$ is computed at $1799$
values of $\alpha\in (0.09, 0.11)$.\\
Figure \ref{figure4} is a graph of $h (\alpha , 10^4, 4\cdot 10^5)$ at $1314$ uniformly
distributed random  values of $\alpha$ in the interval $(0.29, 0.3)$: here 
the non-monotone character of the function $\alpha \mapsto
h(T_\alpha)$ is quite evident. A magnification of Figure \ref{figure5} corresponding to values $\alpha
\in (0.278, 0.281)$, showed in Figure \ref{figure6}, suggests that the same phenomenon
occurs at the end of each of the plateaux exhibited in Figure \ref{figure3}.

\section{Natural extension for $\alpha=\frac{1}{r}$} \label{NE}

In the case $\alpha \in (0,\sqrt{2}-1]$, the structure of the domain $D_\alpha$ of the natural extension for $T_\alpha$ seems to be much more intricate than for $\alpha> \sqrt{2}-1$. Here we find the exact expression for $D_{\alpha}$ and the invariant density of $T_{\alpha}$ when $\alpha \in \left\{\frac{1}{r},r \in \mathbb{N}\right\}$.

\subsection{The by-excess continued fraction map}
\begin{figure}[htb]
\begin{center}
\includegraphics[width=0.6\textwidth, height=0.6\textwidth, angle=270]{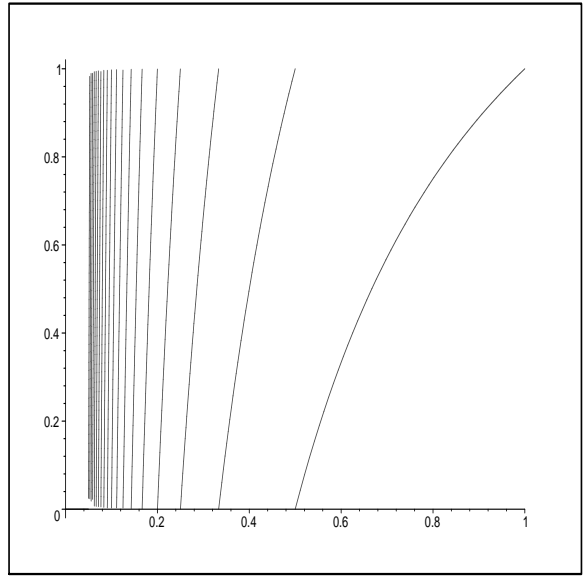}
\caption{Graph of the map $M_0$}
\end{center}
\end{figure}
Before stating our main theorem, we introduce some notations. 
In the following paragraphs we will often refer to the \emph{by-excess} continued fraction expansion of a number, that is the expansion related to the map $M_0(x)=-\frac{1}{x}+\pint{\frac{1}{x}+1}$, $M_0: [0,1] \to [0,1]$. To simplify notations, we will omit the minus signs and use brackets:
\begin{displaymath}
\langle a_0,a_1, a_2, \ldots\rangle \doteqdot \cfrac{1}{a_0- \cfrac{1}{a_1- \cfrac{1}{a_2 - \dotsb}}}, \qquad a_i \in \{2,3,4,\ldots\}
\end{displaymath}		
We will denote a non-integer remainder $x>1$ by a semicolon: 
\begin{displaymath}
\langle a_0,a_1, \ldots,a_n;x\rangle \doteqdot \cfrac{1}{a_0- \cfrac{1}{a_1- \cfrac{}{ \ddots-\cfrac{1}{a_n-\cfrac{1}{x}}}}}, \qquad a_i \in \{2,3,4,\ldots\}
\end{displaymath}	
We also recall that the by-excess expansion of any real number $y \in (0,1)$ is infinite, and that 
\begin{displaymath}
y=\langle a_1,a_2, a_3,\ldots \rangle \in \mathbb{Q} \Rightarrow \exists i \text{ s. t. } \forall j \geq i, \; a_j=2  
\end{displaymath}

\subsection{Reflection rules} \label{translation}
We begin by making some preliminary observations on the relation between the symbolic dynamics of the map $M_0$ and the reflection map $x \mapsto 1-x$ on $[0,1]$, which reveal a sort of \textquotedblleft duality\textquotedblright between the digit $2$ and the digits greater than $2$, and will prove very useful to construct a \textquotedblleft dual\textquotedblright fibred system for $T_\alpha$ in the sense of Schweiger \cite{Sc}.\bigskip 
\par Let $x=\langle a_1,a_2,a_3,\ldots\rangle \in [0,1]$. We would like to determine the by-excess continued fraction expansion of $1-x$. Since the general solution to this problem turns out to be quite complicated, we will only describe a single step of the algorithm, that is, we will suppose to have computed the first $i$ digits of the expansion and the remainder:
\begin{align}
&1-x=\cfrac{1}{{a'_1}- \cfrac{1}{a'_2- \cfrac{1}{\ddots -\cfrac{1}{a'_i -(1-z)}}}}, \qquad z \in [0,1)\notag \\
&z=\langle h_1,h_2,h_3,\ldots\rangle, \qquad h_i\geq 2 \notag
\end{align} 
We want to determine the first digit of the remainder $1-z$.
For reasons that will become clear later, we will treat any sequence of the kind
\begin{displaymath}
\underbrace{2,2,\ldots,2}_{n}
\end{displaymath}
as a single digit. \\
We will make use of the following well-known identity (see for example \cite{MCM}) that can be easily proved by induction on $n$:
\begin{equation} \label{identity}
1-\cfrac{1}{n+\cfrac{1}{y-1}}=\langle\underbrace{2,2,\ldots,2}_{n-1};y\rangle \qquad  \forall y \in \mathbb{R}
\end{equation}  
There are three separate cases to consider: 
\begin{itemize}
\item If $h_1, h_2 \geq 3$, then  from the identity (\ref{identity}) with $n=h_1$ and
\begin{displaymath}
y-1=-\left( h_2-\cfrac{1}{h_3-\ldots}\right), \qquad \frac{1}{y}=-\cfrac{1}{h_2-1-\cfrac{1}{h_3-\ldots}}
\end{displaymath}
we get
\begin{align*}
&1-\langle h_1,h_2,h_3,\ldots\rangle=\big\langle\underbrace{2,\ldots,2}_{h_1-2};2+\langle h_2-1,h_3,\ldots\rangle \big\rangle=\\
&=\big\langle\underbrace{2,\ldots,2}_{h_1-2};3-\big(1-\langle h_2-1,h_3,\ldots\rangle\big) \big\rangle
\end{align*}
We sum up our observations in the following
\medskip
\par \begin{tabular}{|p{0.85\textwidth}|}
\hline
\medskip
\textbf{Rule $1$.\;}
\vspace{-2ex}If $h_1, h_2 \geq 3$,
\begin{equation} 
1-\langle h_1,h_2,h_3,\ldots\rangle=\left\langle\underbrace{2,\ldots,2}_{h_1-2},3;\frac{1}{\left(1-\langle h_2-1,h_3,\ldots\rangle\right)}\right\rangle \notag
\end{equation}
\\
\hline
\end{tabular}
\vspace{2ex}
\item If $z=\langle h_1,\underbrace{2,\ldots,2}_n,h_3,\ldots\rangle$ with $h_1,h_3 \geq 3$, then  
\begin{align*}
1-z=\left\langle \underbrace{2,\ldots,2}_{h_1-2};2+\frac{1}{1-\langle \underbrace{2,\ldots,2}_{n-1},h_3,\ldots\rangle}\right\rangle
\end{align*}
We want to use the identity (\ref{identity}), with
\begin{align*} 
&y-1=-\frac{1}{\langle \underbrace{2,\ldots,2}_{n},h_3,\ldots\rangle}=-\left(2-\langle \underbrace{2,\ldots,2}_{n-1},h_3,\ldots\rangle\right), \\  &-\frac{1}{y}=\frac{1}{1-\langle \underbrace{2,\ldots,2}_{n-1},h_3,\ldots \rangle}
\end{align*}
Observe that 
\begin{align*}
&1-\langle\underbrace{2,\ldots,2}_{n-1},h_3,\ldots\rangle=\frac{1}{n+\langle h_3-1,h_4,\ldots\rangle} \\
& \Rightarrow \frac{1}{1-\langle \underbrace{2,\ldots,2}_{n-1},h_3,\ldots\rangle} =n+1-\left(1-\langle h_3-1,h_4,\ldots\rangle \right)
\end{align*}
In conclusion, we find \medskip
\par \begin{tabular}{|p{0.85\textwidth}|}
\hline
\medskip
\textbf{Rule $2$.\;}
\vspace{-2ex}
If $h_1, h_3 \geq 3$,
\begin{equation} 
1-\langle  h_1,\underbrace{2,\ldots,2}_{n},h_3,h_4,\ldots \rangle=\left\langle \underbrace{2,\ldots,2}_{h_1-2},n+3; \frac{1}{1-\langle h_3-1,h_4,\ldots\rangle}\right\rangle \notag
\end{equation}
\\
\hline
\end{tabular}
\vspace{2ex}
\item If $z=\langle \underbrace{2,\ldots,2}_{n},h_2,\ldots\rangle,\, h_2 \geq 3$, then using again the identity (\ref{identity}) with 
\begin{displaymath}
y=\frac{1}{\langle h_2,h_3,\ldots\rangle}, \qquad \frac{1}{y-1}=\langle h_2-1,h_3,\ldots\rangle
\end{displaymath}
we find 
\begin{align*}
&1-\langle \underbrace{2,\ldots,2}_{n},h_2,h_3,\ldots\rangle=\frac{1}{n+1+\langle h_2-1,h_3,\ldots\rangle}=\\
&=\frac{1}{n+2-\left(1-\langle h_2-1,h_3,\ldots\rangle \right)}
\end{align*}
{}\\
\begin{tabular}{|p{0.9\textwidth}|}
\hline
\vspace{1ex}
\textbf{Rule $3$.\;}
\vspace{-2ex}
If $h_2 \geq 3$,
\begin{equation*}
1-\langle \underbrace{2,\ldots,2}_{n},h_2,h_3,\ldots\rangle=\left \langle n+2;\frac{1}{1-\langle h_2-1,h_3,\ldots\rangle}\right \rangle
\end{equation*}
\\
\hline
\end{tabular}
\end{itemize}\ \\
Notice that we have taken into consideration all the possible cases. Also remark that Rule 1 and Rule 2 guarantee that in the new digits $h_i'$ a sequence of twos is never followed by another. \\
Let $\alpha=\frac{1}{r}$, for a fixed $r \geq 3$. Observe that $T_{\alpha}(\alpha)=0$, and 
\begin{equation*}
T_{\alpha}^i(\alpha-1)= \frac{-(r-i-1)}{r-i}\leq 0 \qquad \text{ for } i=0, \ldots, r-2
\end{equation*} 

Let $\beta$ be the fixed point for $M_0$ corresponding to the branch $r+1$, and $\xi=\frac{1}{r-\beta}$:
\begin{equation} \label{beta}
\begin{split}
&\beta=\frac{r+1-\sqrt{(r+1)^2-4}}{2}=\langle r+1,r+1,r+1, r+1,\ldots\rangle \\  
&\xi= \frac{2}{r-1+ \sqrt{(r+1)^2-4}}=\langle r, r+1, r+1, r+1, \ldots\rangle \\
\end{split} 
\end{equation}

Then 
\begin{equation}
\label{1-beta}
1-\beta=\langle \underbrace{2,\ldots,2}_{r-1},3,\overline{\underbrace{2,\ldots,2}_{r-2},3}\rangle, \qquad   
1-\xi= \langle \overline{\underbrace{2,\ldots,2}_{r-2},3}\rangle 
\end{equation}

\subsection{Domain of the natural extension} \label{definition}
Let $n\geq1$, and define
\begin{multline}
H_n^+=\Bigg \{(h_1,h_2,\ldots,h_n) \Bigg| h_1 \in \{2,(2,2),\ldots,(\underbrace{2,2,\ldots,2}_{r-1})\} \cup \{3,4,\ldots,r\},  \\
 h_2,\ldots,h_{n} \in \{2,(2,2),\ldots,(\underbrace{2,2,\ldots,2}_{r-2})\} \cup \{3,4,\ldots,r,r+1\},\\
 \text{ and such that } h_i=(\underbrace{2,\ldots,2}_s)  \Rightarrow h_{i+1} \geq 3 \Bigg \} \notag
\end{multline}
\begin{multline}
H_n^-=\left \{(h_1,h_2,\ldots,h_{n-1}) \Bigg| h_1, h_2,\ldots,h_{n} \in \right.  \{2,(2,2),\ldots,(\underbrace{2,2,\ldots,2}_{r-2})\} \cup \\
 \cup \{3,4,\ldots,r,r+1\},\; \text{ and such that } h_i=(\underbrace{2,\ldots,2}_s) \Rightarrow h_{i+1} \geq 3 \Bigg \} \notag
\end{multline}
Moreover, for $i=2,3,\ldots,r-1$ define
\begin{multline}
H_n^i=\Bigg \{(h_1,h_2,\ldots,h_n) \Bigg| h_1 \in \{2,(2,2),\ldots,(\underbrace{2,2,\ldots,2}_{r-1-i})\} \cup \{3,4,\ldots,r+1\},  \\
 h_2,\ldots,h_{n} \in \{2,(2,2),\ldots,(\underbrace{2,2,\ldots,2}_{r-2})\} \cup \{3,4,\ldots,r,r+1\}, \\
 \text{ and such that } h_i=(\underbrace{2,\ldots,2}_s) \Rightarrow h_{i+1} \geq 3 \Bigg \} \notag
\end{multline}
Also define
\begin{align*}
&{\hat{H}_n^+}=\left\{ (h_1,h_2,\ldots,h_n) \in H_n^+ \;|\; h_n \geq 3\right \}, \\
&{\hat{H}_n^-}=\left\{ (h_1,h_2,\ldots,h_n) \in H_n^- \;|\; h_n \geq 3\right \}, \\
&{\hat{H}_n^i}=\left\{ (h_1,h_2,\ldots,h_n) \in H_n^i \;|\; h_n \geq 3\right \}, \quad i=2,3,\ldots,r-1 
\end{align*}

Let $V_i(x)=\frac{1}{i-x}$ denote the inverse branches of $M_0$, and
\begin{displaymath}
V_{\big(\underbrace{2,\ldots,2}_s\big)}(x)\doteqdot(\underbrace{V_2 \circ V_2 \circ \cdots \circ V_2}_{s})(x) 
\end{displaymath}
Define
\begin{align*}
&B^+=\bigcup_{n=1}^{\infty} \bigcup_{(h_1,h_2,\ldots,h_n) \in \hat{H}_n^+} (V_{h_1} \circ V_{h_2}\circ \cdots \circ V_{h_n})((1-\xi,1)),
\intertext{and similarly}
&B^-=\bigcup_{n=1}^{\infty} \bigcup_{(h_1,h_2,\ldots,h_n) \in \hat{H}_n^-} (V_{h_1} \circ V_{h_2}\circ \cdots \circ V_{h_n})((1-\xi,1)),\\
&B^i=\bigcup_{n=1}^{\infty} \bigcup_{(h_1,h_2,\ldots,h_n) \in \hat{H}_n^i} (V_{h_1} \circ V_{h_2}\circ \cdots \circ V_{h_n})((1-\xi,1)),\qquad i=2,\ldots,r-1
\end{align*}
Finally, let $E, B, D \subset \mathbb{R}^2$ be defined as follows:
\begin{align*}
& E=\bigcup_{i=1}^{r-1} \left(\left[-\frac{i}{i+1},-\frac{(i-1)}{i}\right] \times \left[0,M_0^{i-1}(1-\xi)\right]\right) \cup \left(\left[0,\frac{1}{r}\right] \times [0,1-\beta]\right), \\
&B=\bigcup_{i=2}^{r-1} \left(\left[-\frac{i}{i+1},-\frac{(i-1)}{i}\right] \times B^i\right)
\cup \left(\left[-\frac{1}{2},0\right] \times B^- \right) \cup \left (\left[0,\frac{1}{r}\right] \times B^+ \right), \\
&D=E \setminus B
\end{align*}
Remark that we have omitted the dependence on $r$ of the sets $B^+, B^-, B^i, E, D$ for simplicity of notation.

\begin{figure}[htb]
\begin{center}
\includegraphics[width=1\textwidth]{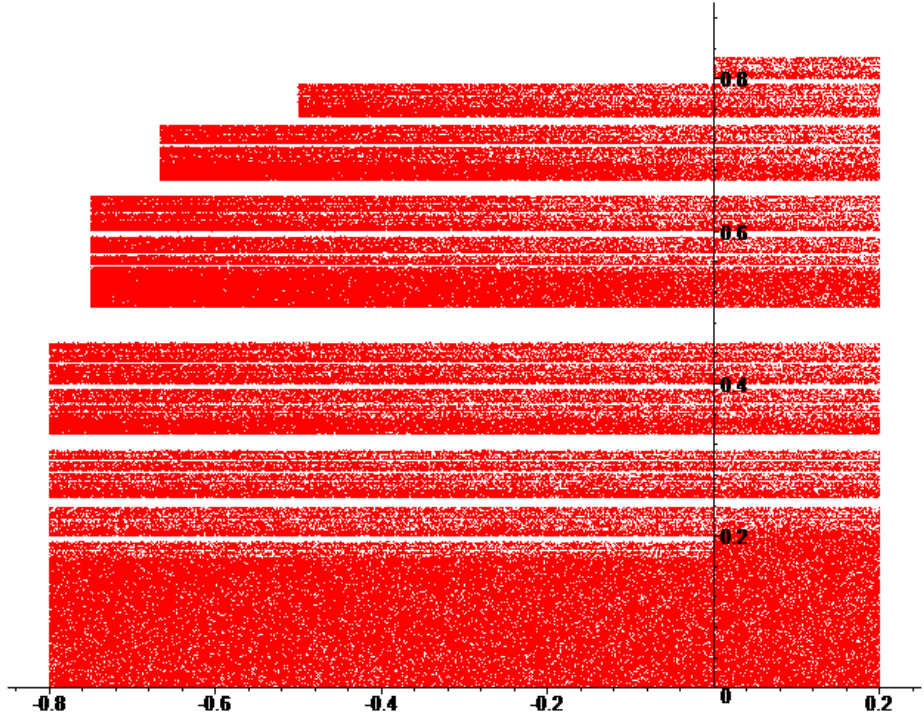}
\caption{A computer simulation for the domain $D$ when $r=5$}
\end{center}
\end{figure}

\begin{theorem} [Natural extension for $\alpha=\frac{1}{r}$] \label{EN}
Let $\alpha=\frac{1}{r}$, $r \geq 3$ be fixed, and let $D \subset \mathbb{R}^2$ be defined as in Paragraph \ref{definition}. Let $k(x)=\pint{\abs{\frac{1}{x}}+1-\alpha}$, and 
\begin{equation} \label{T_bar}
\overline{T}_\alpha(x,y)=\left(T_\alpha(x),\frac{1}{k(x)+\sign(x)y}\right)
\end{equation}
Then $\overline{T}_\alpha: D \to D$ is well defined, one-to-one and onto, and it preserves the density $K_\alpha(x,y)=\frac{1}{C_\alpha} \frac{1}{(xy+1)^2}$, where $C_\alpha=\int_D \frac{1}{(xy+1)^2} dx dy$. In other words, $\overline{T}_\alpha :D \to D$ is a natural extension for $T_\alpha$. 
\end{theorem}

Here the reader should remark that the domain $D$ and the function $k$ also depend on $\alpha$. In the following paragraphs, however, we will write $\overline{T}$ instead of $\overline{T}_\alpha$ for the sake of simplicity.\\  
To prove Theorem \ref{EN} we shall need the following two lemmas:

\begin{lemma} \label{well_defined}
Let $z=\langle h_1,h_2,\ldots,h_n;y \rangle$, where $y>2$ is a real number and $n \geq 1$. Then $1-z$ is of the form $
\left\langle h'_1,h'_2,\ldots,h'_m;\frac{1}{1-1/(y-1)}\right\rangle$, \text{and}
\begin{align*}
&(h_1,h_2,\ldots,h_n) \in H_n^- \Rightarrow (h'_1,h'_2,\ldots,h'_m) \in H_m^+, \\
&(h_1,h_2,\ldots,h_n) \in H_n^+ \Rightarrow (h'_1,h'_2,\ldots,h'_m) \in H_m^-
\end{align*}  
\end{lemma}

\begin{lemma} \label{complementary}
\begin{equation} \label{disjoint}
B^+ \cup \left(1-B^-\right)=\left[\xi,1-\beta\right] \quad (\hspace{-2.3ex}\mod 0)
\end{equation}
and their union is disjoint.
\end{lemma}

\begin{proof}[\textbf{Proof of Lemma \ref{well_defined}}]
Suppose that $(h_1,h_2,\ldots,h_n) \in H_n^-$. From the application of the Rules 1-3, it is straightforward to check that after a suitable number of steps in the algorithm we will obtain a remainder of the form $\frac{1}{1-1/(y-1)}>1$.
We need to verify that at each step of the reflection algorithm described in Paragraph \ref{translation} the newly introduced digits in the by-excess expansion are in accordance with the definition of $H_n^+$. 
We will consider separately the first step and the ensuing ones. 
In the first step, we will have $h_1 \in \{3,\ldots,r,r+1\}$ or $h_1 \in \{2,(2,2),\ldots,(\underbrace{2,2,\ldots,2}_{r-2})\}$.\\ 
If $h_1 \geq 3, h_2 \geq 3$,  applying Rule 1 we get 
\begin{equation} \tag{a}
1-\langle h_1, h_2, h_3,\ldots \rangle = \left \langle \underbrace{2,2,\ldots,2}_{h_1-2},3;\frac{1}{1-\langle h_2-1, h_3, \ldots \rangle} \right\rangle  
\end{equation}
where $h_1-2 \leq r-1$.\\ 
If $h_1 \geq 3, h_2=(\underbrace{2,2,\ldots,2}_{n}),\; n \leq r-2$, using Rule 2 we find
\begin{equation} \tag{b}
1-\langle h_1, h_2, h_3, \ldots \rangle = \left \langle \underbrace{2,\ldots,2}_{h_1-2},n+3; \frac{1}{1- \langle h_3-1, h_4,\ldots \rangle}\right \rangle  
\end{equation}
where $0 \leq h_1-2 \leq r-1, n+3 \leq r+1$.\\
Lastly, for $h_1 \in (\underbrace{2,2,\ldots,2}_{n}),\; n \leq r-2, h_2 \geq 3$, we have
\begin{equation} \tag{c}
1-\langle \underbrace{2,\ldots,2}_{n}, h_2, h_3, \ldots\rangle = \left\langle n+2; \frac{1}{1-\langle h_2-1,h_3,\ldots \rangle} \right\rangle
\end{equation} 
where $n+2 \leq r$ as needed. In all three cases we found an admissible initial segment for $H_m^+$.\\ 
The subsequent steps can be treated in a similar way, although we have to take into account the ways in which the remainder $\langle h_{i+1}, h_{i+2},\ldots \rangle$ from the original sequence has been modified by the reflection rules.
More precisely: if $h_{i+1} \geq 3$ it will be replaced by $h_{i+1}-1 \in \{2,3,\ldots,r\}$; thus when $h_{i+1}-1 \geq 3$, applying Rules 1 and 2, we will find 
\begin{displaymath}
(\underbrace{2,\ldots,2}_{h_{i+1}-3})
\end{displaymath}
as the next digit, with $1 \leq h_{i+1}-3 \leq r-2$, which is admissible for $H_n^+$. Moreover, when $h_{i+1}-1=2$ and $h_{i+2}$ is a sequence of twos they will be considered as a single digit, and it is possible to obtain the sequence
\begin{displaymath}
(\underbrace{2,\ldots,2}_{r-1})
\end{displaymath}
which gives the new digit $r+1$ when we apply Rule 3. We have thus completed the proof for $(h_1,h_2,\ldots,h_n) \in H_n^-$. \\
When switching the roles of $H_n^+$ and $H_n^-$, we can follow the same basic outline. We briefly list the few differences that the reader can easily check for himself: if $(h_1,h_2,\ldots,h_n) \in H_n^+$, 
\begin{itemize}
\item in $(a)$ and $(b)$, we find $h_1 \leq r \Rightarrow h_1-2 \leq r-2$
\item in $(c)$, $n \leq r-1 \Rightarrow n+2 \leq r+1$ 
\end{itemize}
and so the reflected sequence is in accordance with the definition of $H_n^-$. 
\end{proof}

Before moving on to the next Lemma, we make a few observations.\\ 
First of all, notice that since the inverse branches $V_i: x \mapsto \frac{1}{i-x}$ of $M_0$ are all non-decreasing functions, from the by-excess expansions of a sequence of reals we can obtain full knowledge of their ordering. In fact,
\begin{equation} \label{ordering}
\begin{split}
\langle h_1,h_2,\ldots,h_n,\ldots\rangle &< \langle h'_1,h'_2,\ldots,h'_n,\ldots\rangle\\ 
 \Updownarrow \\
\exists i\geq 1 \text{ s. t. } \forall j<i,\; h_j&=h'_j \text{ and } h_i>h'_i
\end{split}
\end{equation}\ \\
Recalling the expansions of $\beta,1-\beta,\xi, 1-\xi$ from equations (\ref{beta}) and (\ref{1-beta}), it follows that $B^- \subset [\beta,1-\xi]$ and $B^+ \subset [\xi,1-\beta]$, and moreover these are the minimal intervals containing $B^+$ and $B^-$: for example, the sequence
\begin{displaymath}
(V_r \circ V_{r+1} \circ V_{r+1} \circ \cdots \circ V_{r+1})(x),\; x \in (1-\xi,1),
\end{displaymath}   
goes arbitrarily close to $\xi$ as the number of pre-images grows. \\
We also observe that if $x \in (1-\xi,1)$, its by-excess expansion must be of the form
\begin{displaymath}
x=\langle \underbrace{2,2,\ldots,2}_{r-1},h_r,h_{r+1},\ldots \rangle, \quad h_r \geq 2
\end{displaymath} 

\begin{proof}[\textbf{Proof of Lemma \ref{complementary}}]
We first want to prove that $B^+$ and $1-B^-$ are disjoint. Let $x \in B^-$; then there exists $l \geq 1$ such that $M_0^l(x) \in (1-\xi,1),\; M_0^j(x) \in [\beta,1-\xi] \; \forall j<l$. Observe that 
\begin{displaymath}
z \in (1-\xi,1) \Rightarrow z=\left \langle \underbrace{\underbrace{2,\ldots,2}_{r-2},3,\ldots,\underbrace{2,\ldots,2}_{r-2},3}_k,\underbrace{2,\ldots,2}_{r-1},\ldots \right \rangle, \quad k\geq 0
\end{displaymath}

Equivalently, for some $i \geq 1$ we have $x= \langle h_1,h_2,\ldots, h_{i},\underbrace{2,\ldots,2}_{n},\ldots \rangle$,
 where
\begin{displaymath}
 n \geq r-1,\;h_{i} \geq 3, \;(h_1,h_2,\ldots,h_{i}) \in \hat{H}_i^-, \; \;(h_1,h_2,\ldots,h_{i-1}) \in H_{i-1}^-%
\end{displaymath}
Then from Lemma \ref{well_defined} we get 
\begin{displaymath}
1-x= \big \langle h'_1,\ldots,h'_m;\frac{1}{1-\langle h_{i}-1,\underbrace{2,\ldots,2}_{n},\ldots \rangle} \big \rangle,
\end{displaymath}
with $(h'_1,\ldots,h'_m) \in H_m^+$,
and applying Rule 2 (or Rule 3 if $h_{i}=3$), we find
\begin{align*}
1-x=\left \langle h'_1,\ldots,h'_m,\underbrace{2,\ldots,2}_{h_{i}-3},n+3;z \right\rangle, \; n+3 \geq r+2,\;\; 0 \leq h_i-3 \leq r-2	
\end{align*}
Observe that $(h'_1,\ldots,h'_m,\underbrace{2,\ldots,2}_{h_{i}-3}) \in H_{m+1}^+$ (or to $H_m^+$ if $h_i-3=0$), but clearly $(h'_1,\ldots,h'_m,\underbrace{2,\ldots,2}_{h_{i}-3},n+3)$ does not belong to $B^+$ because it contains the forbidden digit $n+3$. Since none of the iterates of $1-x$ up to that point belongs to $(1-\xi,1)$, we find that $1-x \notin B^+.$\bigskip 

\par Next we want to show that $B^+ \cup (1-B^-)=(\xi,1-\beta)$.\\
Let $x=\langle h_1,h_2,\ldots \rangle \in (\xi,1-\beta)\setminus B^+$. We must prove that for almost every such $x$ we have $1-x \in B^-$. 
We have to consider two cases:
\begin{enumerate}
\item $\forall n\geq 1, (h_1,\ldots,h_{n}) \in H_n^+$, and so none of the iterates $M_{0}^{n-1}(x)$ belongs to $(1-\xi,1)$ 
\item For some $i$, the by-excess expansion of $x$ contains a forbidden digit $h_i$: either $h_i=(\underbrace{2,\ldots,2}_{n}), n \geq r$ or $h_i \geq r+1$ when $i=1$, or $h_i=(\underbrace{2,\ldots,2}_{n}), n \geq r-1$ or $h_i \geq r+2$ when $i>1$.
\end{enumerate}
However, observe that since the first condition entails in particular that all the elements $h_i$ in the by-excess expansion of $x$ should be bounded, it is satisfied only for a set of Lebesgue measure $0$, and therefore it is negligible for our purposes (equivalently, recall that $M_0$ is ergodic).\\
Next, observe that $x<1-\beta$ implies that the digit $2$ cannot appear $r$ consecutive times in the initial segment of the by-excess expansion of $x$, and $x> \xi$ implies $h_1 \leq r$. Let $i$ be the minimum integer such that $\forall j<i,\; (h_1,\ldots,h_j) \in H_j^+$ and $(h_1,\ldots,h_i) \notin H_i^+$ (we have just seen that $i>1$). Then $h_i$ cannot be of the form $(\underbrace{2,\ldots,2}_{n}), n \geq r-1$ because then $\langle h_i,h_{i+1},\ldots \rangle >1-\xi$ and $x$ would belong to $B^+$.  
The only case left to consider is then $h_i \geq r+2$. Equivalently, one of the iterates $M_{0}^{i-k},\; k \geq 0$ of $x$ is of the form $\langle \underbrace{r+1,\ldots,r+1}_k,r+2,\ldots \rangle < \beta$.
Applying Lemma 1 with $n=i-k-1>1$, $\frac{1}{y}=\big \langle \underbrace{r+1,\ldots,r+1}_k,r+2,\ldots \big \rangle$, we get
$1-x=\left \langle h'_1,\ldots,h'_m;\frac{1}{1-1/(y-1)}\right \rangle, \quad (h'_1,\ldots,h'_m) \in H_m^-$. 
Now observe that $\beta=\frac{1}{r+1-\beta} \Rightarrow \frac{1}{\beta}-1=r-\beta=\frac{1}{\xi}$. 
Then 
$y-1>\frac{1}{\beta}-1>\frac{1}{\xi} \Rightarrow 1-\frac{1}{y-1}>1-\xi$.
Now if $h'_m \geq 3$, we have $(h'_1,\ldots,h'_m) \in \hat{H}_m^-$ and $1-x \in B^-$.
But if $h'_m=(\underbrace{2,\ldots,2}_s)$, we have $h'_{m-1}\geq 3$ and $\left \langle h'_m;\frac{1}{1-1/(y-1)}\right \rangle$ is still greater than $1-\xi$, and again $1-x \in B^-$ (observe that $m>1$, otherwise $1-x=\left \langle h'_m;\frac{1}{1-1/(y-1)}\right \rangle>1-\beta$). 
\end{proof}

\subsection{Proof of Theorem \ref{EN}}
First of all, we observe that Lemma 2 implies that $\overline{T}$ is one-to-one on $D$.
In fact, suppose that $\overline{T}(x_1,y_1)=\overline{T}(x_2,y_2)$. 
Since $y_2 \in [0,1]$, we must have $k(x_2) \in \{k(x_1)-1,k(x_1),k(x_1)+1\}$.
\begin{itemize}
\item If $k(x_1)=k(x_2)$ and $\sign(x_1)=\sign(x_2)$, then obviously $x_1=x_2, y_1=y_2$. 
\item If $k(x_1)=k(x_2)$ and $\sign(x_1)=-\sign(x_2)$, we find $y_1=-y_2$, which is  possible only for $\{y_1=y_2=0\}$, a negligible set. 
\item Lastly, if $x_1>0,x_2<0$ and $k(x_2)=k(x_1)+1$, we get $y_2=1-y_1$. But $(x_1,y_1) \in D \Rightarrow y_1\in [0,\xi]\cup ([\xi,1-\beta]\setminus B^+) \Rightarrow y_2 \in B^- \cup (1-\xi,1) \Rightarrow (x_2,y_2) \notin D$. Thus $\overline{T}$ is one-to-one (mod $0$).\end{itemize} 
\begin{figure}[p]
\begin{center}
\includegraphics[width=0.7\textwidth]{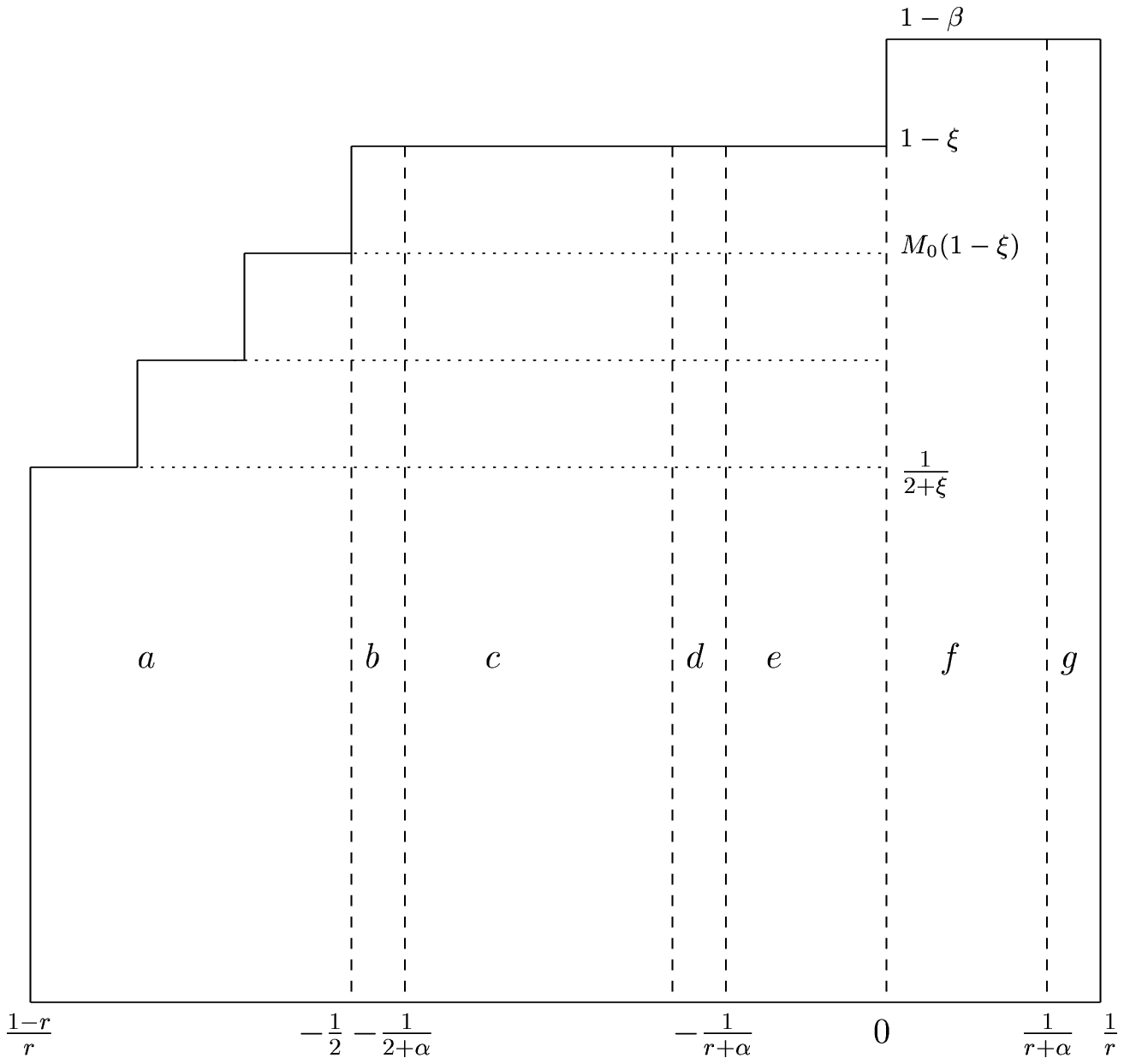}
\caption{A simplified diagram showing the blocks (a)-(g) in the domain}
\end{center}
\end{figure}
\begin{figure}[p]
\begin{center}
\includegraphics[width=0.7\textwidth]{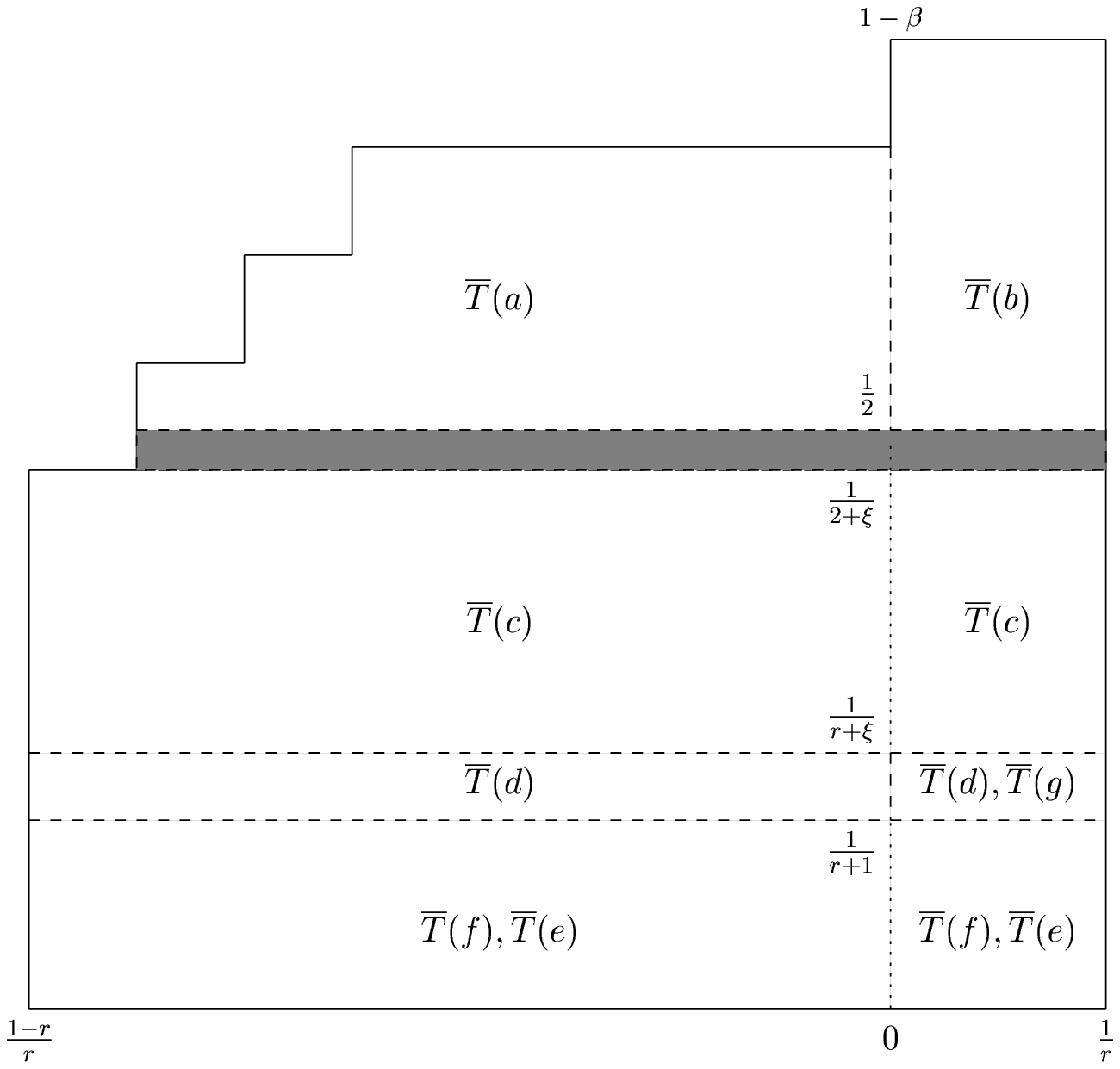}
\caption{A simplified diagram showing the images with respect to $\overline{T}$ of the blocks (a)-(g)}
\end{center}
\end{figure}
Then we can write $\overline{T}(D \setminus B)=\overline{T}(D) \setminus \overline{T} (B)$. Now it is quite straightforward to check that $\overline{T}(D)=D$. In fact, recalling that
\begin{align*}
&\left\{ x>0 \;|\; k(x)=n \right \}=\left(\frac{1}{n+\alpha},\frac{1}{n-1+\alpha}\right],\quad n > r \\
&\left\{ x>0 \;|\; k(x)=r \right \}=\left(\frac{1}{r+\alpha},\alpha\right],\quad n > r \\
& \left\{ x<0 \;|\; k(x)=n \right \}=\left[-\frac{1}{n-1+\alpha},-\frac{1}{n+\alpha}\right), \quad n > 2 \\
& \left\{ x<0 \;|\; k(x)=2 \right \}=\left[\alpha-1,-\frac{1}{2+\alpha}\right)%
\end{align*}
we find 
\begin{multline} 
\overline{T}\left(\left[-\frac{i}{i+1},-\frac{(i-1)}{i}\right] \times \left(\left[0,M_0^{i-1}(1-\xi)\right]\setminus B^i\right)\right)= \\ =\left[-\frac{(i-1)}{i},-\frac{(i-2)}{i-1}\right] \times \left(\left[\frac{1}{2},M_0^{i-2}(1-\xi)\right]\setminus B^{i-1}\right),\, i=2,\ldots,r-1, \tag{a} 
\end{multline}
\begin{multline}
\overline{T}\left(\left[-\frac{1}{2},-\frac{1}{2+\alpha}\right] \times \left( \left[0,1-\xi\right]\setminus B^- \right)\right)=[0,\alpha]\times \left(\left[\frac{1}{2},1-\beta\right]\setminus B^+\right), \tag{b} 
\end{multline}
\begin{multline}
\overline{T}\left(\left[-\frac{1}{n-1+\alpha},-\frac{1}{n+\alpha}\right] \times \left([0,1-\xi]\setminus B^-\right)\right)=\\
=[\alpha-1,\alpha] \times \left(\left[\frac{1}{n},\frac{1}{n-1+\xi}\right]\setminus B^-\right),\quad n=3,\ldots,r \tag{c}
\end{multline}
Here we observe that $B^+\cup \left[\frac{1}{r},\frac{1}{2+\xi}\right]=B^-\cup \left[\frac{1}{r},\frac{1}{2+\xi}\right]=B^i \cup \left[\frac{1}{r},\frac{1}{2+\xi}\right]$ for $i=2,3,\ldots,r-1$. Also remark that the rectangles $\left[-\frac{(r-2)}{r-1},\alpha\right]\times \left[\frac{1}{2+\xi},\frac{1}{2}\right]$ and $[\alpha-1,\alpha] \times \left[\frac{1}{n+\xi},\frac{1}{n}\right]$ for $n=3,\ldots,r$ both belong to $B$.   
\begin{multline} \label{d} \tag{d}
\overline{T}\left(\left[-\frac{1}{r+\alpha},-\frac{1}{r+1+\alpha}\right] \times \left([0,1-\xi]\setminus B^-\right)\right)=\\=\left([\alpha-1,0] \times \left(\left[\frac{1}{r+1},\frac{1}{r+\xi}\right]\cap D \right)\right)  \cup \\ \cup   \left([0,\alpha] \times \left(\left[\frac{1}{r+1},\frac{1}{r+\xi}\right]\setminus V_{r+1}(B^-)\right)\right)
\end{multline} 
(Here we wanted to highlight the fact that $B^+ \cap \left[\frac{1}{r+1},\frac{1}{r+\xi}\right]=\varnothing$.) 
\begin{multline}\tag{e}
\overline{T}\left(\left[-\frac{1}{n-1+\alpha},-\frac{1}{n+\alpha}\right] \times \left([0,1-\xi]\setminus B^-\right)\right)=\\=[\alpha-1,\alpha] \times \left(\left[\frac{1}{n},\frac{1}{n-1+\xi}\right]\setminus V_{n}(B^-)\right),\quad n
\geq r+2%
\end{multline}
\begin{multline}\tag{f}
\overline{T}\left(\left[\frac{1}{n-1+\alpha},\frac{1}{n-2+\alpha}\right] \times \left([0,1-\beta]\setminus B^+\right)\right)=\\=[\alpha-1,\alpha] \times \left(\left[\frac{1}{n-\beta},\frac{1}{n-1}\right]\setminus V_{n-1}^+(B^+)\right), \quad n\geq r+1 
\end{multline}
where we set $V_{n}^+(x)= \frac{1}{n+x}$. 
\begin{multline}\tag{g}
\overline{T}\left(\left[\frac{1}{r+\alpha},\alpha\right]\times \left([0,1-\beta]\setminus B^+\right) \right)=\\=[0,\alpha] \times \left(\left[\frac{1}{r+1-\beta},\frac{1}{r}\right]\setminus V_r^+(B^+)\right)
\end{multline}
To conclude the proof observe that 
\begin{displaymath}
B^+ \cup (1-B^-)=[\xi,1-\beta] \Rightarrow V_r^+(B^+)=\left[\frac{1}{r+1-\beta},\frac{1}{r+\xi}\right]\setminus V_{r+1}(B^-)
\end{displaymath}
which together with (\ref{d}) proves that $\overline{T}$ is onto.\\ 
Finally, the fact that $K(x,y)$ is invariant for $\overline{T}$ can be easily checked through the change of variables formula: we observe that the determinant of the Jacobian for $\overline{T}$ is respectively $\frac{1}{x^2(k(x)+y)^2}$ when $x>0$ and $\frac{1}{x^2(k(x)-y)^2}$ when $x<0$. Then for any $A \subset D$, if we put $A^+=A\cap \{x>0\}$ and $A^-=A \cap \{x<0\}$, we have
\begin{multline*}
\int_{\overline{T}(A)} K(x,y) dx dy=
\frac{1}{C_\alpha} \left( \int_{A^+} \frac{1}{u^2(1/u+v)^2} du dv + \right.\\
\left. +\int_{A^-} \frac{1}{u^2(-1/u-v)^2} du dv \right)= \frac{1}{C_\alpha}\left(
\int_{A^+} \frac{du dv}{(1+uv)^2}+ \int_{A^-} \frac{du dv}{(1+uv)^2} \right) \qedhere
\end{multline*} 

\subsection{Invariant densities and entropy for $\alpha=\frac{1}{r}$}
Since $\pi_1 \circ \overline{T}_{\alpha}=T_{\alpha} \circ \pi_1$, where $\pi_1$ is simply the projection on the first coordinate, the invariant density for $T_\alpha$ is obtained by integrating $K_\alpha(x,y)$ with respect to the second coordinate. 
Given a sequence $(h_1,h_2,\ldots,h_n)$, define
\begin{align*}
&a(h_1,h_2,\ldots,h_n)=\frac{1}{\left \langle h_1,h_2,\ldots,h_n;\frac{1}{1-\xi} \right \rangle}>1\\
&b(h_1,h_2,\ldots,h_n)=\frac{1}{\left \langle h_1,h_2,\ldots,h_n-1 \right \rangle}>1 %
\end{align*} %
Let 
\begin{align*}
&\psi_+(x)=\sum_{n=1}^{\infty} \sum_{(h_1,\ldots,h_n) \in H^+} \left(\frac{1}{x+b(h_1,\ldots,h_n)}-\frac{1}{x+a(h_1,\ldots,h_n)} \right), \\
&\psi_-(x)=\sum_{n=1}^{\infty} \sum_{(h_1,\ldots,h_n) \in H^-} \left(\frac{1}{x+b(h_1,\ldots,h_n)}-\frac{1}{x+a(h_1,\ldots,h_n)} \right), \\
&\psi_i(x)=\sum_{n=1}^{\infty} \sum_{(h_1,\ldots,h_n) \in H^i} \left(\frac{1}{x+b(h_1,\ldots,h_n)}-\frac{1}{x+a(h_1,\ldots,h_n)} \right)
\end{align*}
for $i=2,\ldots,r-1$, and observe that
\begin{displaymath}
\int_{a}^{b} \frac{1}{(1+xy)^2} dy= \frac{1}{x+\frac{1}{b}}-\frac{1}{x+\frac{1}{a}}
\end{displaymath}
It follows that, for a suitable normalization constant $c_\alpha$, 
\begin{align*}
&\rho_{\alpha}=\frac{\psi_\alpha(x)}{c_\alpha}=\frac{1}{c_\alpha}\left(\sum_{i=2}^{r-1} \left( \chi_{\left[-\frac{i}{i+1},-\frac{(i-1)}{i}\right]}(x) \left( \frac{1}{x+\frac{1}{M_0^{i-1}(1-\xi)}}-\psi_i(x)\right)\right)+ \right. \\ &\left. +\chi_{\left[-\frac{1}{2},0\right]}(x) \left(\frac{1}{x+\frac{1}{1-\xi}}-\psi_-(x)\right)+
\chi_{\left[0,\frac{1}{r}\right]}(x) \left(\frac{1}{x+\frac{1}{1-\beta}}-\psi_+(x)\right)\right)
\end{align*}
is an invariant density for $T_\alpha$. 

\begin{remark} 
Even in the case $\alpha=\frac{1}{r}$ the domain of the natural extension seems too complicated to allow for a direct computation of the entropy. However, as far as Corollary \ref{zero} is concerned, it is pro\-ba\-bly possible to prove a much stronger result. In fact Nakada \cite{Nak} showed that in the case of $\frac{1}{2} \leq \alpha \leq 1$, the integral $-2\int_{\alpha-1}^{\alpha} \log \abs{x} \int_{D_\alpha(x)} K(x,y) dy dx$ (where $D_{\alpha}(x)$ are the vertical sections of the domain of the natural extension) is constant and equal to $\frac{\pi^2}{6}$. We conjecture that the same should be true for $0<\alpha<\frac{1}{2}$. 
\end{remark}

\section{Appendix} \label{Appendix}
In this section we collect the proofs of some lemmas that were used in the main part of the article.

\subsection{Bounded Distortion} \label{bounded_distortion}
\begin{proof}[\textbf{Proof of Proposition \ref{bound_dist}}]
Observe that $\exists k>0$ such that $\forall I^{\varepsilon}_{j} \in \mathcal{P}, \forall x, y \in I^{\varepsilon}_{j}$,
\begin{equation*}
\abs{\frac{T_{\alpha}'(x)}{T_{\alpha}'(y)}-1} \leq k \abs{T_{\alpha}(x)-T_{\alpha}(y)}
\end{equation*}
In fact, if $x,y \in I^{\varepsilon}_{j,\alpha}$, then 
\begin{align*}
&\abs{\frac{T_{\alpha}'(x)}{T_{\alpha}'(y)} -1} \frac{1}{\abs{T_{\alpha}(x)-T_{\alpha}(y)}}= \abs{\frac{y^2}{x^2}-1} \frac{\abs{xy}}{\abs{x-y}} \leq \abs{\frac{y}{x}}\abs{x+y}\leq 4
\end{align*}
Let $n \geq 1,\; I_{\eta}^{(n)} \in \mathcal{P}^{(n)}$,\;$x,y \in I_{\eta}^{(n)}$. Define $\lambda=\sup \abs{\frac{1}{T_{\alpha}'}}=(1-\alpha)^2$: then 
\begin{multline} \label{bd}
\log \abs{\frac{(T_{\alpha}^n)'(y)}{(T_{\alpha}^n)'(x)}}=\sum_{i=0}^{n-1} \log \abs{\frac{T_{\alpha}'(T_{\alpha}^i(y))}{T_{\alpha}'(T_{\alpha}^i(x)}} \leq \sum_{i=0}^{n-1} \abs{\frac{T_{\alpha}'(T_{\alpha}^i(y))}{T_{\alpha}'(T_{\alpha}^i(x))}-1} \leq \\
\leq 4 \sum_{i=0}^{n-1} \abs{T_{\alpha}^{i+1}(y)-T_{\alpha}^{i+1}(x)}=4 \sum_{i=1}^n \abs{T_{\alpha}^{i}(y)-T_{\alpha}^{i}(x)} \leq \\
 \leq 4 \sum_{i=1}^n \lambda^{n-i} \abs{T_{\alpha}^{n}(y)-T_{\alpha}^{n}(x)} \leq 
4 \sum_{i=0}^\infty \lambda^i = \frac{4}{1-(1-\alpha)^2}=C_2 
\end{multline}
Then $\abs{\frac{(T_{\alpha}^n)'(y)}{(T_{\alpha}^n)'(x)}} \leq e^{C_2}= C_1$. Let $I_{\eta}^{(n)}$ be a full cylinder: $T_{\alpha}^n(I_{\eta}^{(n)})=I_{\alpha}$.  
Now consider any measurable set $B$: 
\begin{multline} \label{ratio}
\frac{m(B)}{m(I_{\alpha})}=\frac{\int\limits_{V_\eta(B)}\abs{(T_{\alpha}^n)'(y)}dy}{\int\limits_{I_{\eta}^{(n)}}\abs{(T_{\alpha}^n)'(x)}dx} \leq \frac{m(V_\eta(B)) \sup\limits_{y \in I_{\eta}^{(n)}} \abs{(T_{\alpha}^n)'(y)}}{m(I_{\eta}^{(n)}) \inf\limits_{x \in I_{\eta}^{(n)}} \abs{(T_{\alpha}^n)'(x)}} \leq C_1 \frac{m(V_\eta(B))}{m(I_{\eta}^{(n)})} \\
\Rightarrow m(V_{\eta}(B)) \geq m(B) \frac{m(I_{\eta}^{(n)})}{C_1}, 
\end{multline}
which concludes the proof. 	
\end{proof}

\subsection{Exactness} \label{exact}
To prove the exactness of the $T_\alpha$, $\alpha \in \left(0,\frac{1}{2}\right)$, we follow the same argument of H. Nakada (\cite{Nak}, Theorem 2).
The crucial property we need in order to prove Lemma \ref{exactness} is 

\begin{proposition}\label{Borel}
The family of the cylinder sets $I_{\eta}^{(n)}=(\omega_1,\ldots,\omega_n) \in \mathcal{P}^{(n)}$ such that $T_\alpha^n(I_{\eta}^{(n)})=I_\alpha$ generates the Borel sets.
\end{proposition}

(Here we write $\omega_i=(j_i,\varepsilon_i)$ for brevity). 

\begin{proof}[\textbf{Proof of Proposition \ref{Borel}}]
Consider the sets
\begin{displaymath}
E_n=\{ (\omega_1,\ldots,\omega_n) \;|\; T_\alpha(\omega_1) \neq I_\alpha, T_{\alpha}^{2}(\omega_1,\omega_2)\neq I_\alpha,\ldots,T_{\alpha}^n(\omega_1,\ldots,\omega_n) \neq I_\alpha\}
\end{displaymath}
and let $M_n=m\left(\bigcup_{I_{\eta}^{(n)} \in E_n} I_{\eta}^{(n)}\right)$. Consider the orbits of the endpoints with respect to $T_\alpha$:
\begin{align*}
\alpha=(a_1,a_2,a_3,\ldots), \qquad \alpha-1=(b_1,b_2,b_3,\ldots)
\end{align*}
Then $E_1=\{(a_1),(b_1)\}$, and 
\begin{multline*}
E_n=\{(\omega_1,\ldots,\omega_n) \in \mathcal{P}^{(n)} \;|\; (\omega_2,\ldots,\omega_n) \in E_{n-1} \text { and } \omega_1=a_1 \text{ or } b_1\} \\ \cup  \{(a_1,a_2,\ldots,a_n),(b_1,b_2,\ldots,b_n)\}
\end{multline*}
In fact if $\omega_1 \notin \{a_1,b_1\}$, we would have $T_\alpha(\omega_1) =I_\alpha$; moreover, if $(\omega_2,\ldots,\omega_n) \neq (a_2,\ldots,a_n)$, the monotonicity of $T_\alpha$ on $(a_1)$ implies that either $(\omega_2,\ldots,\omega_n) \cap T_\alpha(a_1) = \varnothing$, or $(\omega_2,\ldots,\omega_n) \subseteq T_{\alpha}(a_1)$. In this last case $T_{\alpha}^n(a_1,\omega_2,\ldots,\omega_n)=T_{\alpha}^{n-1}(\omega_2,\ldots,\omega_n)$.
So we get
\begin{align*}
M_n \leq ((1-\alpha)^2+\alpha^2)M_{n-1}+m((a_1,\ldots,a_n) \cup (b_1,\ldots,b_n)),
\end{align*}
and since $(1-\alpha)^2+\alpha^2<1$ and $m(w_1,\ldots,w_n)$ vanishes as $n \to \infty$, we have $M_n \to 0$ as $n \to \infty$, that is, $E=\{x \;|\; \forall n,\; T_{\alpha}^n(I_{\eta}^{(n)}(x)) \neq I_\alpha\}$
has Lebesgue measure $0$, where  $I_{\eta}^{(n)}(x)$ is the cylinder in $\mathcal{P}^{(n)}$ containing $x$. Then, recalling that $T_{\alpha}$ is non-singular, $m(T_{\alpha}^{-n}(E))$ is also $0$ for all $n \geq 0$, and so $m\left(\bigcup_{n} T_{\alpha}^{-n}(E)\right)=0$. That is, for almost all $x$ there is a subsequence $\{n_i\}$ such that $T^{n_i}(I_{\eta}^{n_i}(x))=I_{\alpha}$ for all $i \in \mathbb{N}$. Then for almost all $x$, $\forall U$ open neighborhood of $x$ we can find $n$ and a full cylinder $x \in I_{\eta}^{(n)} \subset U$.   
\end{proof}

\begin{proof}[\textbf{Proof of Lemma \ref{exactness}}]
We have just proved that the full cylinders generate the Borel sets. Then a sufficient and necessary condition for exactness, due to Rohlin \cite{Ro}, is the following: $\exists C>0$ such that $\forall n, \; \forall I_{\eta}^{(n)}$ full cylinder of rank $n$, $\forall X \subset I_{\eta}^{(n)}$,
\begin{equation} \label{Rohlin_c}
\mu_{\alpha}(T_{\alpha}^n(X)) \leq C \frac{\mu_{\alpha}(X)}{\mu_{\alpha}(I_{\eta}^{(n)})}
\end{equation}
We recall that the $T_{\alpha}$ satisfy the \emph{bounded distortion} property: 
Then, recalling that the density of $\mu_{\alpha}$ with respect to the Lebesgue measure is bounded from above and from below by constants, we get for some constant $C$,
\begin{equation*}
\mu_{\alpha}(V_\eta(B)) \geq \frac{1}{C} \mu_{\alpha}(B)\mu_{\alpha}(I_{\eta}^{(n)}),
\end{equation*}
that is, Rohlin's characterization (\ref{Rohlin_c}).
\end{proof}

\subsection{Continuity} \label{app}
\begin{proof}[\textbf{Proof of Lemma \ref{continuity}}]
Since $\sup \abs{\tilde{\rho}_{\alpha_n}} \leq K$, $\Var \tilde{\rho}_{\alpha_n} \leq K$ $\forall n$, we can apply the following theorem:

\famoustheo[Helly's Theorem]
{Let $\{\rho_n\}$ be a sequence in $BV(I)$ such that:
\begin{enumerate}
\item $\sup \abs{\rho_n} \leq K_1 \quad \forall n$, 
\item $\Var \rho_n \leq K_2 \quad \forall n$
\end{enumerate}   
Then there exists a subsequence $\rho_{n_k}$ and a function $\rho \in BV(I)$ such that
$\rho_{n_k} \xrightarrow{\;L^1\;} \rho$, $\rho_{n_k} \to \rho$ almost everywhere, and
\begin{align*}
& \sup \abs{\rho} \leq K_1, \qquad \Var \rho \leq K_2
\end{align*}}

Thus we can find a subsequence $\{\tilde{\rho}_{\alpha_{n_k}}\}$ converging in the $L^1$ norm and almost everywhere to some function $\rho_\infty$ such that $\sup \abs{\rho_\infty} \leq K$, $\Var \rho_\infty \leq K$. We want to show that $\rho_\infty=\tilde{\rho}_{\bar{\alpha}}$: we observe that it is sufficient to show that $\rho_\infty$ is an invariant density for $A_{\bar{\alpha}}=A_{\bar{\alpha},\bar{\alpha}}=T_{\bar{\alpha}}$, and then use the uniqueness of the invariant density.   
To simplify notations, we will write $\alpha_k$ for $\alpha_{n_k}$, $\tilde{\rho}_k$ for $\tilde{\rho}_{\alpha_{n_k}}$, and $A_k$ for $A_{\alpha_{n_k},\bar{\alpha}}$. \\
Our goal is to show that $\forall B \subseteq I_{\bar{\alpha}}$,\;$\int \chi_B(A_{\bar{\alpha}}(x)) \rho_\infty(x) dx=\int \chi_B(x) \rho_\infty(x) dx$.
Observe that every $\chi_B(x)$ belongs to $L^1(I_{\bar{\alpha}})$ and so can be approximated arbitrarily well by compactly supported $C^1$ functions with respect to the $L^1$ norm. Then it will be sufficient to prove that $\forall \varphi \in C^1$ with compact support contained in $I_{\bar{\alpha}}$,
\begin{equation} \label{tesi}
\abs{\int \varphi\left(A_{\bar{\alpha}}(x)\right) \rho_\infty(x) dx - \int \varphi(x) \rho_\infty(x) dx} =0   
\end{equation}
Observe that $\left |\int \varphi(A_{\bar{\alpha}}(x))\rho_{\infty}(x) dx - \int \varphi (x) \rho_{\infty}(x) dx  \right| \leq I_1+I_2+I_3$, with $I_1,I_2,I_3$ given below:
\begin{align*}
&I_1=\left |\int \varphi(A_{\bar{\alpha}}(x))\rho_{\infty}(x) dx - \int \varphi (A_{\bar{\alpha}}(x)) \tilde{\rho}_k(x) dx  \right| \leq \norm{\varphi}_{\infty} \norm{\tilde{\rho}_k-\rho_\infty}_{L^1} \\
&I_3=\left |\int \varphi(A_k(x))\tilde{\rho}_k(x) dx - \int \varphi (x) \rho_\infty(x) dx  \right| = \\ 
&=\left |\int \varphi(x)\tilde{\rho}_k(x) dx - \int \varphi (x) \rho_\infty(x) dx  \right| \leq
\norm{\varphi}_\infty \norm{\tilde{\rho}_k-\rho_\infty}_{L^1}  
\end{align*}
which vanish as $k \to \infty$. Finally, $I_2=\int\abs{\varphi(A_{\bar{\alpha}}(x))- \varphi(A_k(x))}\tilde{\rho}_k(x)dx$
is bounded by $K \int \abs{\varphi(A_k(x))-\varphi(A_{\bar{\alpha}}(x))} dx
$, and we need to show that 
\begin{equation} \label{caso2}
\int \abs{\varphi(A_k(x))-\varphi(A_{\bar{\alpha}}(x))} dx \to 0 \text{ when } k \to \infty
\end{equation}
Recall that for $x \in J^{-}_{j,\alpha_k}=\left[\frac{-1}{j-1+\alpha_k}+\bar{\alpha}-\alpha_k,-\frac{1}{j+\alpha_k}+\bar{\alpha}-\alpha_k\right)$, $A_k(x)=-\frac{1}{x-\bar{\alpha}+\alpha_k}$ $-j+\bar{\alpha}-\alpha_{k}$, and for $x \in J^{-}_{j,\bar{\alpha}} =\left[-\frac{1}{j-1+\alpha},-\frac{1}{j+\alpha}\right)$, $A_{\bar{\alpha}}(x)=-\frac{1}{x}-j$.
We will examine in detail the case $\alpha_{k}< \bar{\alpha} \; \forall k$\;,\;$x<\bar{\alpha}-{\alpha}_{k}$; the other cases can be dealt with in a similar way. In this case, $0<\frac{1}{j+\alpha_k}-\frac{1}{j+\bar{\alpha}}<\bar{\alpha}-\alpha_k$, and if $j <\frac{1}{\sqrt{\bar{\alpha}-{\alpha}_{k}}}=N(k)$, then  $-\frac{1}{j-1+\alpha_k}+\frac{1}{j+\bar{\alpha}}<\frac{\alpha_k -\bar{\alpha}-1}{j^2}<\alpha_k-\bar{\alpha}$ and so 
\begin{equation} \label{N}
-\frac{1}{j-1+{\alpha}_{k}}+\bar{\alpha}-{\alpha}_{k} < -\frac{1}{j+\bar{\alpha}} < -\frac{1}{j+{\alpha}_{k}}+\bar{\alpha}-{\alpha}_{k}
\end{equation}
$I_{N(k)}=\bigcup_{j\geq N(k)} J^{-}_{j,\alpha_k}$ contains the set in which condition (\ref{N}) isn't satisfied, and its measure $m(I_{N(k)})=\sum_{j=N(k)}^{\infty} \abs{\frac{1}{(j-1+{\alpha}_{k})(j+{\alpha}_{k})}}$ vanishes when $k \to \infty$. 
Given $\varepsilon''>0$, choose $\bar{k}$ such that $m(I_{N(\bar{k})})< \varepsilon''$, and let $k\geq \bar{k}$. Define 
\begin{align*}
&\xi^{-}_{j}= \left[-\frac{1}{j-1+\alpha_k}+\bar{\alpha}-\alpha_{k},-\frac{1}{j+\bar{\alpha}}\right), \; 
\eta^{-}_{j}= \left[-\frac{1}{j+\bar{\alpha}},-\frac{1}{j+\alpha_k}+\bar{\alpha}-\alpha_k\right)\\
&\text{when } 2<j\leq N(\bar{k}), \text{ and } \quad \xi^{-}_2=\left[\bar{\alpha}-1,-\frac{1}{2+\alpha_{k}}\right)
\end{align*}
Then we can split the integral (\ref{caso2}) in three parts in the following way: 
\begin{multline}
\int_{\bar{\alpha}-1}^{\bar{\alpha}} \abs{\varphi(A_k (x))-\varphi(A_{\bar{\alpha}}(x))} dx \leq
\int\limits_{I_{N(\bar{k})}} \abs{\varphi(A_k (x))-\varphi(A_{\bar{\alpha}}(x))} dx+ \\
+\sum_{j=3}^{N(\bar{k})} \left( \int_{\eta^{-}_{j}} \abs{\varphi(A_k (x))-\varphi(A_{\bar{\alpha}}(x))} dx \right) + \sum_{j=2}^{N(\bar{k})} \left(
 \int_{\xi^{-}_j} \abs{\varphi(A_k (x))-\varphi(A_{\bar{\alpha}}(x))} dx\right) \notag 
\end{multline}
The first integral in this expression is bounded by $2 \varepsilon'' \norm{\varphi}_\infty$. Moreover, the measures of the sets $\eta^{-}_{j}$ tend uniformly to $0$ when $k \to \infty$: 
\begin{align*}
&m(\eta^{-}_{j}) \leq \abs{\bar{\alpha}-\alpha_{k}}+\frac{\abs{\bar{\alpha}-\alpha_k}}{(j+\bar{\alpha})(j+\alpha_{k})} \leq C_1 (\bar{\alpha}-\alpha_k) \\
&\Rightarrow \sum_{j=3}^{N(\bar{k})} \int_{\eta^{-}_j} \abs{\varphi(A_k (x))-\varphi(A_{\bar{\alpha}}(x))} dx \leq N(\bar{k}) 2\norm{\varphi}_\infty C_1 (\bar{\alpha}-\alpha_k)
\end{align*}
Finally, $m(\xi^{-}_j)\leq \frac{C_2}{j^2}+\abs{\bar{\alpha}-\alpha_k}\leq \frac{C_3}{j^2}$ when $k \geq \bar{k},\; j<N(\bar{k})$, and for $x \in \xi^{-}_{j}$, $x \leq -\frac{1}{j+\bar{\alpha}}$ and $x-\bar{\alpha}+\alpha_k <-\frac{1}{j+\alpha_k}$, therefore 
\begin{multline*}
\abs{A_k(x)-A_{\bar{\alpha}}(x)} =\abs{-\frac{1}{x}+\frac{1}{x-\bar{\alpha}+\alpha_{k}}-\bar{\alpha}+\alpha_{k}} \leq\\ \leq \abs{\bar{\alpha}-\alpha_{k}}+\frac{\abs{\bar{\alpha}-\alpha_{k}}}{\abs{x(x-\bar{\alpha}+\alpha_k)}}\leq \abs{\bar{\alpha}-\alpha_k}(1+(j+1)^2)
\end{multline*}
Since $\varphi$ is $C^1$ on a compact interval, it is also lipschitzian for some Lipschitz constant $L_\varphi$, and 
\begin{align*}
&\sum_{j=2}^{N(\bar{k})}\int_{\xi^{-}_j} \abs{\varphi(A_k (x))-\varphi(A_{\bar{\alpha}}(x))} dx \leq
\sum_{j=2}^{N(\bar{k})} m(\xi^{-}_j) L_\varphi \abs{A_k(x)-A_{\bar{\alpha}}(x)} \leq \\
& \leq \sum_{j=2}^{N(\bar{k})} \frac{C_3(1+(j+1)^2)}{j^2}L_\varphi \abs{\bar{\alpha}-\alpha_{k}}\leq C_4 N(\bar{k}) \abs{\bar{\alpha}-\alpha_{k}} \leq C_4 \sqrt{\abs{\bar{\alpha}-\alpha_{k}}} 
\end{align*}
when $k$ is large. This establishes the claim that the third integral vanishes when $x<\bar{\alpha}-\alpha_{k}$. In the case $x>\bar{\alpha}-\alpha_{k}$ we have similar estimates: for $j<\frac{1}{\sqrt{\abs{\bar{\alpha}-\alpha_{k}}}}$, we have 
\begin{equation*}
\frac{1}{j+\bar{\alpha}}<\frac{1}{j+\alpha_k}+\bar{\alpha}-\alpha_{k}<\frac{1}{j-1+\bar{\alpha}}
\end{equation*}
and we can define the intervals
\begin{equation*}
\gamma^{+}_{j}=\left(\frac{1}{j+\alpha_k}+\bar{\alpha}-\alpha_{k},\frac{1}{j-1+\bar{\alpha}}\right],\; \delta^{+}_{j}=\left(\frac{1}{j-1+\bar{\alpha}},\frac{1}{j-1+\alpha_{k}}+\bar{\alpha}-\alpha_{k}\right]
\end{equation*}
We have $m(\delta^{+}_{j})\leq C_5 \abs{\bar{\alpha}-\alpha_k}$, $m(\gamma^{+}_{j})\leq \frac{C_5}{j^2}$, and
\begin{equation*}
 \abs{A_k(x)-A_{\bar{\alpha}}(x)} \leq C_7 j^2 \abs{\alpha_k-\bar{\alpha}} \text{ for } x \in \gamma^{+}_{j}
\end{equation*}
Finally, we leave it to the reader to check that the case $\bar{\alpha} < \alpha_k$ can be treated in same way.  
Thus we can conclude that (\ref{caso2}) holds. \\ 
Therefore we have shown that $\rho_\infty=\tilde{\rho}_{\bar{\alpha}}$. This is also true if we extract a converging sub-subsequence from any subsequence of $\tilde{\rho}_{\alpha_n}$, and so $\tilde{\rho}_{\alpha_n} \to \tilde{\rho}_{\bar{\alpha}}$ both in $L^1$ and almost everywhere for $n \to \infty$. This completes the proof of Lemma \ref{continuity}. 
\end{proof}

\end{document}